\documentclass[a4paper]{article}

\usepackage{amsmath,amsthm,amssymb}
\usepackage{graphicx}
\usepackage{epstopdf}
\usepackage{listings}
\usepackage{hyperref}

\usepackage{algorithm,algpseudocodex}

\usepackage{listings}

\DeclareMathOperator{\diag}{diag}

\theoremstyle{remark}

\newcommand{\A}{\mathbf A}
\newcommand{\B}{\mathbf B}
\newcommand{\C}{\mathbf C}
\newcommand{\CC}{\mathbb C}
\newcommand{\D}{\mathbf D}
\newcommand{\dd}{\mathbf d}
\newcommand{\hatDta}{\hat{\mathbf D}_t^\alpha}
\newcommand{\Dta}{\mathbf D_t^\alpha}
\newcommand{\Dx}{\mathbf D_x}
\newcommand{\E}{\mathbf E}
\newcommand{\e}{\mathbf e}
\newcommand{\Eta}{\mathbf E_t^\alpha}
\newcommand{\hatEta}{\hat{\mathbf E}_t^\alpha}
\newcommand{\f}{\mathbf f}
\newcommand{\hatf}{\hat{\mathbf f}}
\newcommand{\F}{\mathbf F}

\newcommand{\G}{\mathbf G}
\newcommand{\HH}{\mathbf H}

\newcommand{\I}{\mathbf I}
\newcommand{\U}{\mathbf U}
\newcommand{\Uin}{\mathbf U_{inner}}
\newcommand{\M}{\mathbf M}

\newcommand{\R}{\mathbb{R}}
\newcommand{\N}{\mathbb{N}}

\newcommand{\x}{\mathbf x}

\newcommand{\dig}{\mathtt{dig}}
\newcommand{\X}{\mathbf X}

\newcommand{\rl}{I_t^{\alpha}}

\newcommand{\cpt}{D_t^{\alpha}}

\newcommand{\err}{\operatorname{err}}

\newcommand{\errhatDta}{\err_{\hatDta}}
\newcommand{\errDta}{\err_{\Dta}}

\newcommand{\mx}[1]{\|#1\|_{\max}}

\begin{document}
	
\title{Numerical approximation of Caputo-type advection-diffusion equations in one and multiple spatial dimensions via shifted Chebyshev polynomials}

\author{Francisco de la Hoz\thanks{Department of Mathematics, Faculty of Science and Technology, University of the Basque Country UPV/EHU, Barrio Sarriena S/N, 48940 Leioa, Spain.} \ and Peru Muniain\thanks{Department of Applied Mathematics, Escuela de Ingeniería de Bilbao, University of the Basque Country UPV/EHU, Plaza Ingeniero Torres Quevedo, 1, 48013 Bilbao, Spain. Corresponding author (peru.muniain@ehu.eus).}}

\date{University of the Basque Country UPV/EHU}
	
\maketitle	
	
\begin{abstract} 

In this paper, using a pseudospectral approach, we develop operational matrices based on the shifted Chebyshev polynomials to approximate numerically Caputo fractional derivatives and Riemann-Liouville fractional integrals. In order to make the generation of these matrices stable, we use variable precision arithmetic. Then, we apply the Caputo differentiation matrices to solve numerically Caputo-type advection-diffusion equations in one and multiple spatial dimensions, which involves transforming the discretization of the concerning equation into a Sylvester (tensor) equation. We provide complete Matlab codes, whose implementation is carefully explained. The numerical experiments involving highly oscillatory functions in time confirm the effectiveness of this approach.

\medskip


\end{abstract}

\section{Introduction}

Fractional calculus (see, e.g., \cite{podlubny} for a complete overview) generalizes the ideas of differentiation and integration to noninteger orders. In this paper, among the possible definitions of fractional derivatives, we consider the Caputo fractional derivative, which is defined as follows. Let $\alpha\in(0,\infty)$, such that $0 \le n - 1 < \alpha < n$, for $n\in\mathbb N$. Then, the Caputo fractional derivative of order $\alpha$ of a function $f(t)\in\mathcal C^{n}(0, \infty)$ is given by
\begin{equation}
\label{e:cptf}
\cpt f(t) = \frac{1}{\Gamma(n - \alpha)}\int_0^t\frac{f^{(n)}(\tau)}{(t - \tau)^{1 - n + \alpha}}d\tau.
\end{equation}
Unlike the integer-order derivatives $f'(t)$, $f''(t)$, etc., which are local, fractional derivatives are nonlocal. However, it is immediate to see that
\begin{align}
\label{e:liman}
\lim_{\alpha\to n^-}\cpt f(t) & = \ f^{(n)}(t),
	\\
\nonumber \lim_{\alpha\to (n-1)^+}\cpt f(t) & = \ f^{(n-1)}(t) - f^{(n-1)}(0) \Longrightarrow \lim_{\alpha\to n^+}\cpt f(t) = \ f^{n}(t) - f^{n}(0),
\end{align}
where we have integrated by parts in the first equality. In particular, $\lim_{\alpha\to0^+}\cpt f(t) = f(t) - f(0)$. Therefore, the operator $\cpt$ is not continous at $\alpha\in\mathbb N$.

Likewise, let $\alpha\in(0,\infty)$. Then, the Riemann-Liouville fractional integral of a function $f(t)\in\mathcal C^0(0, \infty)$ is given by
\begin{equation}
\label{e:rlf}
\rl f(t) = \frac{1}{\Gamma(\alpha)}\int_0^t\frac{f(\tau)}{(t - \tau)^{1-\alpha}}d\tau.
\end{equation}
Note that $\rl f(t)$ is not defined at $\alpha = 0$, but, integrating \eqref{e:rlf} by parts, it follows immediately that $\lim_{\alpha\to0^+}\rl f(t) = f(t)$.

In this paper, we are mainly interested in approximating numerically \eqref{e:cptf} in $t\in[0, T]$. However, since the numerical treatment of \eqref{e:rlf} is very similar, we also consider it. Let us mention that there are other possible notations to denote \eqref{e:cptf} and \eqref{e:rlf}. For instance, it is possible to add respectively $C$ (which stands for Caputo) and $RL$ (which stands for Riemann-Liouville) as superscripts or subscripts (by typing, e.g., $^C\cpt f(t)$ and $^{RL}\rl f(t)$). However, since there is no risk of confusion, we simply write $\cpt f(t)$ and $\rl f(t)$.

Among the existing methods (see, e.g., the review \cite{caili2020} and its references), we will focus on the methods based on polynomial interpolation. This can be done, e.g., by means of Lagrange interpolating polynomials, as in \cite{highorder1,highorder2,highorder3}, where $\cpt f(t)$ is approximated in the context of solving numerically Caputo-type advection-diffusion equations, or more recently in \cite{delahozmuniain2024}, where an FFT-based convolution algorithm applied to a modification of the methods in \cite{highorder1,highorder2,highorder3} allows to approximate numerically $\cpt f(t)$ for an extremely large number of nodes in a very efficient way. Moreover, in \cite{delahozmuniain2024}, after discretizing the time and space variables by means of operational matrices, Caputo-type advection-diffusion equations are transformed into a Sylvester equation of the form $\A\cdot\X+\X\cdot\B=\C$, which enables to consider simultaneously the solution all the time instants and all the spatial nodes.

On the other hand, it is also possible to perform the interpolation by means of a spectral method and, in this regard, the so-called shifted Chebyshev polynomials $T_k^*(t)$, $t\in[0,1]$, which will be used in this paper, seem a pretty natural option. These polynomials are just the first-kind Chebyshev polynomials $T_k(t) = \cos(k\arccos(t))$ (see, e.g., \cite{trefethen,boyd} for a complete exposition of their properties) evaluated at $2t - 1$, i.e.,
$$
T_k^*(t) = T_k(2t - 1) = \cos(\arccos(2t - 1)).
$$
Then, given $N\in\mathbb N\cup\{0\}$, we approximate $f(t)$ as
\begin{equation}
	\label{e:ftTkstar}
	f(t) \approx \sum_{k=0}^N\hat f_kT_k^*\left(\frac tT\right), \quad t\in[0, T],
\end{equation}
for a certain $T > 0$, and apply \eqref{e:cptf} and \eqref{e:rlf} to \eqref{e:ftTkstar}, getting respectively
\begin{align}
\label{e:cptfsum}
	\cpt f(t) & \approx \sum_{k=0}^N\hat f_k \cpt \left(T_k^*\left(\frac tT\right)\right), \quad t\in[0, T],
	\\
\label{e:rlfsum}
	\rl f(t) & \approx \sum_{k=0}^N\hat f_k \rl\left(T_k^*\left(\frac tT\right)\right), \quad t\in[0, T].
\end{align}
In what regards the actual implementation of \eqref{e:cptfsum} and/or \eqref{e:rlfsum} and differential equations involving them, it is possible to do it by constructing operational matrices, as in, e.g, \cite{graef2014,wangachena2020,farhood2023}, but even if such matrices are not explictly defined, as in e.g., \cite{hashemi2023}, it is necessary to solve linear systems of equations.

The main inconvenient of using shifted Chebyshev polynomials for the numerical approximation of fractional derivatives and/or integrals is that the resulting methods are in general unstable, except for small values of $N$, which limits their applicability. Moreover, the existing theoretical results are rather involved; for instance, in \cite[Th. 1]{waleed2021}, the authors proved the following:
$$
\cpt T_k(t) \approx \sum_{p = 0}^{M}\Delta_{k,p}^{(\mu)}T_p(t), \quad M \gg 1,
$$
where
\begin{equation*}
\Delta_{k,p}^{(\mu)} = \frac{2 k (-1)^{k + \lceil\alpha\rceil}\Gamma(k + \lceil\alpha\rceil) \Gamma\left(\lceil\alpha\rceil - \alpha + \frac{1}{2}\right)\theta_p}{\Gamma(k - \lceil\alpha\rceil + 1)}{}_4\tilde F_3\left(\genfrac{}{}{0pt}{}{1, \lceil\alpha\rceil-k, k+\lceil\alpha\rceil, \lceil\alpha\rceil - \alpha +\frac{1}{2}}{\lceil\alpha\rceil+\frac{1}{2}, \lceil\alpha\rceil-\alpha-p+1, \lceil\alpha\rceil-\alpha+p+1} \Bigg| 1 \right),
\end{equation*}
with ${}_4\tilde F_3$ denoting the regularized ${}_4F_3$ function, and
$$
\theta_p =
\begin{cases}
\frac12, & p = 0,
\\
1, & p \ge 1.
\end{cases}
$$
Bearing in mind the previous arguments, the main aim of this work is to generate fractional differentiation matrices and fractional integration matrices that are stable for any $N$, for which we will use arbitrary precision arithmetic (vpa) in the intermediate steps. Let us remark that the use of vpa has been used successfully in \cite{cayama2020}, to evaluate numerically the Gaussian hypergeometric function ${}_2F_1$ in the context of the fractional Laplacian, which is a related operator, but, to the best of our knowledge, it has not been used previously to generate fractional operational matrices associated to the shifted Chebyshev polynomials.

The structure of this paper is as follows. In Section~\ref{s:pseudospectral}, we approximate a given function $f(t)$ in terms of shifted Chebyshev polynomails, as in \eqref{e:ftTkstar}, and express \eqref{e:cptfsum} and \eqref{e:rlfsum} by means of operational matrices. In Section~\ref{s:matlab}, we explain carefully how to implement in Matlab \cite{matlab} the matrices developed in Section~\ref{s:pseudospectral}, by means of vpa. In this regard, special care is given to choosing the correct number of digits. Note that, once that the matrices have been obtained by means of vpa, they can be safely cast to 64-bit precision, but, except for small values of $N$, it is not possible to compute correctly the matrices without vpa. In Section \eqref{s:numerical}, we perform several numerical tests, considering different values of $\alpha$ and $N$, obtaining very satisfactory results, specially for $\alpha\in(0,1)$. Finally, in Section \eqref{s:advdif}, as a practical application, following the ideas in \cite{delahozmuniain2024} and \cite{cuestadelahoz2024}, we solve Caputo-type advection-diffusion equations in one and multiple spatial dimensions.

To help understanding the implementation and reproduction of the numerical experiments, complete Matlab codes are offered.

All the codes have been run in a Mainstream A+ Server AS-2024S-TR SUPERMICRO AMD, with 56 cores, 112 threads, 2.75 GHz, and 128 GB of RAM.

\section{A matrix-based pseudospectral method method to approximate numerically $\cpt f(t)$ and $\rl f(t)$}

\label{s:pseudospectral}

Given a regular enough function $f(t)$ defined on $[0, T]$, with $T > 0$, in order to approximate numerically \eqref{e:cptf} and \eqref{e:rlf} by means of \eqref{e:cptfsum} and \eqref{e:rlfsum}, respectively, we consider a pseudospectral method, where, in order to determine the coefficients $\{\hat f_k\}$ in \eqref{e:ftTkstar}, we must impose the equality in \eqref{e:ftTkstar} at $N + 1$ given nodes.

\subsection{Obtention of $\{\hat f_k\}$}

In order to obtain $\{\hat f_k\}$ in \eqref{e:ftTkstar}, we impose the equality at the $N + 1$ extreme points of $T_k^*(t/T)$:
\begin{equation}
\label{e:tj}
t_j = \frac T2\left[1+\cos\left(\frac{j\pi}{N}\right)\right], \quad 0 \le j \le N.
\end{equation}
Therefore, denoting $f_j$ the numerical approximation of $f(x_j)$:
\begin{align}
\label{e:fj}
f_j \equiv \sum_{k=0}^N\hat f_kT_k^*\left(\frac{t_j}{T}\right) = \sum_{k=0}^N\hat f_kT_k\left(2\frac{t_j}{T}-1\right) =  \sum_{k=0}^N\hat f_kT_k\left(\cos\left(\frac{j\pi}{N}\right)\right)
 =  \sum_{k=0}^N\hat f_k\cos\left(\frac{jk\pi}{N}\right), \quad 0 \le j \le N,
\end{align}
i.e., we have a discrete cosine transform. Then, representing the cosines in exponential form:
\begin{align}
\label{e:fxj}
f_j & = \sum_{k=0}^N\hat f_k\frac{e^{ijk\pi/N} + e^{-ijk\pi/N}}2 = \frac12\sum_{k=0}^N\hat f_ke^{2ijk\pi/(2N)} + \frac12\sum_{k=-N}^0\hat f_{-k}e^{2ijk\pi/(2N)}
	\cr
& = \hat f_0 + \frac12\sum_{k=1}^{N-1}\hat f_ke^{2ijk\pi/(2N)} + (-1)^{jk}\hat f_N + \frac12\sum_{k=N+1}^{2N-1}\hat f_{2N-k}e^{2ijk\pi/(2N)}, \quad 0 \le j \le N.
\end{align}
On the other hand, if we replace $j$ by $2N - j$ in the last equation,
\begin{align*}
f_{2N-j} & =  \hat f_0 + \frac12\sum_{k=1}^{N-1}\hat f_ke^{2i(2N-j)k\pi/(2N)} + (-1)^{(2N-j)k}\hat f_N + \frac12\sum_{k=N+1}^{2N-1}\hat f_{2N-k}e^{2i(2N-j)k\pi/(2N)}
	\cr
& =  \hat f_0 + \frac12\sum_{k=N+1}^{2N-1}\hat f_{2N-k}e^{-2ij(2N-k)\pi/(2N)} + (-1)^{jk}\hat f_N + \frac12\sum_{k=1}^{N-1}\hat f_{k}e^{-2ij(2N-k)\pi/(2N)} = f_j,
\end{align*}
where we have substituted $k$ by $2N-k$ in the sums of the second line. Therefore, defining
\begin{equation}
\label{e:gj}
g_j =
\begin{cases}
f_j, & 0 \le j \le N,
	\cr
f_{2N-j}, & N + 1 \le j \le 2N-1,
\end{cases}
\end{equation}
and
\begin{equation}
\label{e:hatgk}
\hat g_k =
\begin{cases}
\hat f_0, & k = 0,
	\cr
\frac12\hat f_k, & 1 \le k \le N - 1,
	\cr
\hat f_N, & k = N,
	\cr
\frac12\hat f_{2N-k}, & N \le k \le 2N - 1,
\end{cases}
\end{equation}
we conclude that
\begin{equation}
\label{e:gjhatgk}
g_j = \sum_{k=0}^{2N-1}\hat g_ke^{2ijk\pi/(2N)} \Longleftrightarrow \hat g_k = \frac{1}{2N}\sum_{j=0}^{2N-1}g_je^{-2ijk\pi/(2N)},
\end{equation}
i.e., $\{g_j\}$ and $\{\hat g_k\}$ are obtained from each other by means of an inverse discrete Fourier transform (IDFT) and a discrete Fourier transform (DFT), respectively. Therefore, to get $\{\hat f_k\}$ from $\{f_j\}$, we define $\{g_j\}$ according to \eqref{e:gj}, perform a DFT on them, which yields $\{\hat g_k\}$, and apply \eqref{e:hatgk}, to get $\{\hat f_k\}$. Likewise, to obtain $\{f_j\}$ from $\{\hat f_k\}$, we define $\{\hat g_k\}$ according to \eqref{e:hatgk}, perform an IDFT on them, which yields $\{g_j\}$, and apply \eqref{e:gj}, to get $\{f_j\}$, but keeping only $j\in\{0, 1, \ldots, N\}$.

\subsubsection{Matrix-based obtention of $\{\hat f_k\}$ from $\{f_j\}$}

Bearing in mind the previous arguments, it is pretty straightforward to create a matrix $\M = [m_{ij}]\in\mathbb R^{(N+1)\times(N+1)}$, such that
\begin{equation}
\label{e:M}
\hatf = \M \cdot \f
\Longleftrightarrow
\f =\M^{-1}\cdot \hatf,
\end{equation}
where
\begin{equation}
\label{e:fhatf}
\f = 
\begin{pmatrix}
	f_0 \\ f_1 \\ \vdots \\ f_N
\end{pmatrix},
	\qquad
\hatf = 
\begin{pmatrix}
	\hat f_0 \\ \hat f_1 \\ \vdots \\ \hat f_N
\end{pmatrix}.
\end{equation}
Indeed, combining \eqref{e:gj} and \eqref{e:gjhatgk},
\begin{equation*}
\hat g_k = \frac{f_{0}}{2N} + \frac{1}{N}\sum_{j=1}^{N-1}f_{j}\cos\left(\frac{jk\pi}{N}\right) + \frac{(-1)^kf_{N}}{2N},
\end{equation*}
and hence, from \eqref{e:hatgk},
\begin{equation*}
\hat f_k =
\left\{
\begin{aligned}
& \frac{f_{0}}{2N} + \frac{1}{N}\sum_{j=1}^{N-1}f_{j} + \frac{f_{N}}{2N}, & & k = 0,
	\cr
& \frac{f_{0}}{N} + \frac{2}{N}\sum_{j=1}^{N-1}f_{j}\cos\left(\frac{jk\pi}{N}\right) + \frac{(-1)^kf_{N}}{N}, & & 1 \le k \le N - 1,
	\cr
& \frac{f_{0}}{2N} + \frac{1}{N}\sum_{j=1}^{N-1}(-1)^jf_j + \frac{1}{2N}f_{N}(-1)^N, && k = N.
\end{aligned}
\right.
\end{equation*}
Therefore, the entries of the symmetric matrix $\M$ are given by
\begin{equation}
\label{e:mjk}
m_{jk} = m_{kj} = \frac{2}{N}\varepsilon_{jk}\cos\left(\frac{jk\pi}{N}\right), \text{ where }
\varepsilon_{jk} = 
\begin{cases}
\frac14, & j,k \in \{0, N\},
	\\
\frac12, & j \in \{0, N\} \text{ and } 1 \le k \le N - 1,
	\\
\frac12, & 1 \le j \le N - 1 \text{ and } k \in \{0, N\},
	\\
1, & 1 \le i,j \le N - 1.
\end{cases}
\end{equation}
Note that the inverse of $\M$ follows trivially from \eqref{e:fj}:
$$
[\M^{-1}]_{jk} = \cos\left(\frac{jk\pi}{N}\right), \quad 0 \le j, k \le N.
$$

\subsection{Construction of the fractional differentiation and integration matrices}

Coming back to the shifted Chebyshev polynomials $\{T_k^*(t)\}$, we are interested in obtaining explicitly its coefficients, i.e., find $\{a_{kl}\}$, such that
\begin{equation}
\label{e:Tkstarakl}
T_k^*(t) = \sum_{l = 0}^{k}a_{kl}t^l.
\end{equation}
Even if there are different ways to do this, we opt for using a recurrence to obtain $\{T_k^*(t)\}$ recursively, in a similar way as with the standard Chebyshev polynomials:
\begin{equation}
\label{e:Chebrecurs}
\begin{cases}
T_0(t) = 1,
	\\
T_1(t) = t,
	\\
T_{k+1}(t) = 2tT_k(t) - T_{k-1}(t), & k \in \mathbb N.
\end{cases}
\end{equation}
Hence, replacing $t$ by $2t - 1$ in \eqref{e:Chebrecurs} and bearing in mind that $T_k^*(t) \equiv T_k(2t-1)$, we have
\begin{equation}
	\label{e:shiftChebrecurs}
	\begin{cases}
		T_0^*(t) = 1,
		\\
		T_1^*(t) = 2t-1,
		\\
		T_{k+1}^*(t) = (4t-2)T_k(t) - T_{k-1}^*(t), & k \in \mathbb N.
	\end{cases}
\end{equation}
In this way, we obtain the values of the coefficients $\{a_{kl}\}$ of $\{T_0^*, T_1^*, \ldots, T_N^*\}$ as defined in \eqref{e:Tkstarakl}, and store them in an upper triangular matrix $\C = [c_{ij}]\in\mathbb Z^{(N+1)\times(N+1)}$, in such a way that $c_{l+1,k+1} \equiv a_{kl}$, i.e., the coefficients of the $k$th polynomial $T_k^*(t)$ are stored at the $(k+1)$th column, and the constant coefficient $a_{k0}$ of $T_k^*(t)$ occupies the first position of that column; the coefficient $a_{k1}$ multiplying $x$, the second position; the coeficient $a_{k2}$ multiplying $x^2$, the third position, and so on. For instance, if $N = 5$, then
\begin{align*}
T_0(t) & = 1,
	\\
T_1(t) & = 2t-1,
	\\
T_2(t) & = 8t^2 - 8t + 1,
	\\
T_3(t) & = 32t^3 - 48t^2 + 18t - 1,
	\\
T_4(t) & =  128t^4 - 256t^3 + 160t^2 - 32t + 1,
	\\
T_5(t) & = 512t^5 - 1280t^4 + 1120t^3 - 400t^2 + 50t - 1,
\end{align*}
from which we get
\begin{equation*}
\C =
\begin{pmatrix}
	a_{00} & a_{10} & a_{20} & a_{30} & a_{40} & a_{50}
	\\
	0 & a_{11} & a_{21} & a_{31} & a_{41} & a_{51}
	\\
	0 & 0 & a_{22} & a_{32} & a_{42} & a_{52}
	\\
	0 & 0 & 0 & a_{33} & a_{43} & a_{53}
	\\
	0 & 0 & 0 & 0 & a_{44} & a_{54}
	\\
	0 & 0 & 0 & 0 & 0 & a_{55}
\end{pmatrix}
=
\begin{pmatrix}
1 & -1 & 1 & -1 & 1 & - 1
	\\
0 & 2 & -8 & 18 & -32 & 50
	\\
0 & 0 & 8 & -48 & 160 & -400
	\\
0 & 0 & 0 & 32 & -256 & 1120
	\\
0 & 0 & 0 & 0 & 128 & -1280
	\\
0 & 0 & 0 & 0 & 0 & 512
\end{pmatrix},
\end{equation*}
note that the structure of the first row comes from the fact that $c_{1,k+1} = a_{k0} = T_k^*(1) = T_k(1) = (-1)^k$, for $0 \le k \le N$.

From a practical point of view, in order to generate $\C$ by using \eqref{e:shiftChebrecurs}, we observe that, if we multiply by $x$ the polynomial corresponding to a given column, this is equivalent to shifting the column one position downwards. Therefore, we create first an all-zero matrix $\C$ of order $N + 1$, and the first two columns, corresponding respectively to the first two equations in \eqref{e:shiftChebrecurs}, are completely determined by assigning $c_{11} = 1$, $c_{12} = -1$, $c_{22} = 2$. Then, the $j$th column, for $3 \le j \le N + 1$, is generated from the $j$th and $(j-1)$th  according to the third equation in \eqref{e:shiftChebrecurs}, by the following recurrence:
\begin{equation}
\label{e:cij}
\begin{cases}
c_{1j} = (-1)^{j-1},
	\cr
c_{ij} = 4c_{i - 1, j - 1} - 2c_{i, j - 1} - c_{i, j - 2}, & 3 \le i \le j \le N + 1.
\end{cases}
\end{equation}
After obtaining $\C$, and hence determining the coefficients $\{a_{kl}\}$ in \eqref{e:Tkstarakl}, we introduce \eqref{e:Tkstarakl} in \eqref{e:ftTkstar}, getting
\begin{equation}
\label{e:ftsumsum}
f(t) \approx \sum_{k=0}^N\sum_{l = 0}^{k}a_{kl}\hat f_k\left(\frac tT\right)^l = \sum_{k=0}^N\sum_{l = 0}^{k}c_{l+1,k+1}\hat f_k\left(\frac tT\right)^l = \sum_{k=0}^N\sum_{l = 0}^{N}c_{l+1,k+1}\hat f_k\left(\frac tT\right)^l,
\end{equation}
where we have used in the last equality that $\C$ is upper triangular. Therefore, in order to apply \eqref{e:cptf} and \eqref{e:rlf} to \eqref{e:ftsumsum}, we consider first their action on $f(t) = (t / T)^l$. In the case of \eqref{e:cptf}, we recall that, when $n, l\in\mathbb N\cup\{0\}$,
$$
\frac{d^n(t^l)}{dt^n} = 
\left\{
\begin{aligned}
& l(l-1)\ldots(l + 1  - n)t^{l-n} = \frac{\Gamma(l + 1)}{\Gamma(l + 1 - n)}t^{l-n}, & & n \le l,
	\cr
& 0, & & n > l.
\end{aligned}
\right.
$$
Hence, introducing $f(t) = (t / T)^l$, with $l\in\mathbb N\cup\{0\}$, in \eqref{e:cptf}, it follows that, when $\alpha > l$, $\cpt (t/T)^l = 0$, and, when $\alpha \le l$,
\begin{align}
\label{e:cpttl}
\cpt \left(\frac tT\right)^l & = \frac{\Gamma(l + 1)}{T^l\Gamma(n-\alpha)\Gamma(l + 1 - n)}\int_0^t\frac{\tau^{l-n}}{(t - \tau)^{1 - n + \alpha}}d\tau
	\cr
& = \frac{\Gamma(l + 1)t^{l-\alpha}}{T^l\Gamma(n-\alpha)\Gamma(l  + 1 - n)}\int_0^1z^{l + 1 - n - 1}(1 - z)^{n - 1 - \alpha}d\tau
	\cr
& = \frac{\Gamma(l + 1)t^{l-\alpha}}{T^l\Gamma(n-\alpha)\Gamma(l  + 1 - n)}B(l + 1 - n, n - \alpha) = \frac{\Gamma(l + 1)}{T^l\Gamma(l + 1 - \alpha)}t^{l-\alpha},
\end{align}
where we have performed the change of variable $\tau = tz$, and $\Gamma(\cdot)$ and $B(\cdot,\cdot)$ denote Euler’s gamma and beta functions, respectively. Note that \eqref{e:cpttl} is formally valid, when $l + 1 - \alpha\not\in-\mathbb Z\cup\{0\}$. In particular, if $\alpha \le l$ is a nonnegative integer,
$$
\cpt \left(\frac tT\right)^l = \frac{l(l - 1)\ldots(l + 1 - \alpha)}{T^l}t^{l-\alpha},
$$
which is precisely the standard, integer-order derivative of order $\alpha$ of $f(t) = (t / T)^l$.

In what regards \eqref{e:rlf}, introducing $f(t) = (t / T)^l$, with $l\in\mathbb N\cup\{0\}$, in \eqref{e:cptf},
\begin{align}
\label{e:rltl}
\rl \left(\frac tT\right)^l & = \frac{1}{T^l\Gamma(\alpha)}\int_0^t\frac{\tau^l}{(t - \tau)^{1-\alpha}}d\tau = \frac{t^{l + \alpha}}{T^l\Gamma(\alpha)}\int_0^1z^{l + 1 - 1}(1 - z)^{-1 + \alpha}dz
	\cr
& = \frac{t^{l + \alpha}}{T^l\Gamma(\alpha)}B(l + 1, \alpha) = \frac{\Gamma(l + 1)}{T^l\Gamma(l + 1 + \alpha)}t^{l + \alpha};
\end{align}
note that \eqref{e:rltl} is formally valid, when $l + 1 + \alpha\not\in-\mathbb Z\cup\{0\}$. Moreover, \eqref{e:rltl} is precisely \eqref{e:cpttl}, after replacing $-\alpha$ by $\alpha$.

Bearing in mind the previous arguments, we approximate numerically \eqref{e:cptf} and \eqref{e:rlf} by applying respectively the operators $\cpt$ and $\rl$ to \eqref{e:ftsumsum}, and evaluating the resulting expressions at $t = t_j$:
\begin{equation}
\label{e:cptftj}
\cpt f(t_j) \approx \sum_{k=0}^N\sum_{l = 0}^{N}c_{l+1,k+1}\hat f_k\cpt\left(\frac{t_j}T\right)^l = \sum_{k=0}^N\left[\sum_{l = \lceil\alpha\rceil}^{N}t_j^{l - \alpha}\frac{\Gamma(l + 1)}{T^l\Gamma(l + 1 - \alpha)}c_{l+1,k+1}\right]\hat f_k, \quad 0 \le j \le N,
\end{equation}
where $\lceil\alpha\rceil \le N$, and
\begin{equation}
\label{e:rlftj}
\rl f(t_j) \approx \sum_{k=0}^N\sum_{l = 0}^{N}c_{l+1,k+1}\hat f_k\rl\left(\frac{t_j}T\right)^l = \sum_{k=0}^N\left[\sum_{l = 0}^{N}t_j^{l + \alpha}\frac{\Gamma(l + 1)}{T^l\Gamma(l + 1 + \alpha)}c_{l+1,k+1}\right]\hat f_k, \quad 0 \le j \le N.
\end{equation}
Observe that, if $p(x)$ is a polynomial of degree lower than or equal to $N$, then \eqref{e:ftsumsum} is exact for $f(x) = p(x)$, and so are \eqref{e:cptftj} and \eqref{e:rlftj}. Moreover, even if the operators $\cpt$ and $\rl$ are defined in \eqref{e:cptf} and \eqref{e:rlf} for positive noninteger values of $\alpha$, both approximations \eqref{e:cptftj} and \eqref{e:rlftj} are well defined for $\alpha\in\N\cup\{0\}$. More precisely, when $\alpha = 0$, $\cpt f(t_j) = \rl f(t_j)$, whereas, when $\alpha\in\N\cup\{0\}$, $\cpt f(t_j)\approx f^{(\alpha)}(t_j)$. The latter is coherent with \eqref{e:liman}, and is due to the fact that we are taking $\lceil\alpha\rceil$ as the lower limit of the inner sum in \eqref{e:cptftj}.

Coming back to \eqref{e:cptftj} and \eqref{e:rlftj}, in order to express them in matrix form, we define, respectively, $\hatDta = [d_{ij}]\in\mathbb C^{(N+1)\times(N+1)}$ and $\hatEta = [e_{ij}]\in\mathbb C^{(N+1)\times(N+1)}$ as
\begin{align}
\label{e:dj1k1}
d_{j+1,k+1} & = \sum_{l = \lceil\alpha\rceil}^{N}t_j^{l - \alpha}\frac{\Gamma(l + 1)}{T^l\Gamma(l + 1 - \alpha)}c_{l+1,k+1}, \quad 0 \le j, k \le N,
	\\
\label{e:ej1k1}
e_{j+1,k+1} & = \sum_{l = 0}^{N}t_j^{l + \alpha}\frac{\Gamma(l + 1)}{T^l\Gamma(l + 1 + \alpha)}c_{l+1,k+1}, \quad 0 \le j, k \le N,
\end{align}
where, in order not to burden the notation, we omit the subscript $_t$ and the superscript $^\alpha$ in $d_{j+1,k+1}$ and $e_{j+1,k+1}$. Then,
\begin{equation}
\label{e:DatEatfj}
\begin{pmatrix}
\cpt f(t_0)
	\\
\cpt f(t_1)
	\\
	\vdots
	\\
\cpt f(t_N)
\end{pmatrix}
\approx
\hatDta \cdot
\hatf,
	\qquad
\begin{pmatrix}
	\rl f(t_0)
	\\
	\rl f(t_1)
	\\
	\vdots
	\\
	\rl f(t_N)
\end{pmatrix}
\approx
\hatEta \cdot
\hatf,
\end{equation}
where $\hatf$ is defined in \eqref{e:fhatf}. Note that it is possible to approximate $\{\cpt f(t_j)\}$ and $\{\rl f(t_j)\}$ directly from $\{f_j\}$ by means of $\M$ defined in \eqref{e:mjk}. Indeed, from \eqref{e:M}:
\begin{equation}
	\label{e:tildeDatEatfj}
	\begin{pmatrix}
		\cpt f(t_0)
		\\
		\cpt f(t_1)
		\\
		\vdots
		\\
		\cpt f(t_N)
	\end{pmatrix}
	\approx
	\Dta \cdot\f,
	\qquad
	\begin{pmatrix}
		\rl f(t_0)
		\\
		\rl f(t_1)
		\\
		\vdots
		\\
		\rl f(t_N)
	\end{pmatrix}
	\approx
	\Eta \cdot \f,
\end{equation}
where $\f$ is defined in \eqref{e:fhatf}, and
\begin{equation}
\label{e:tildeDattildeEat}
\Dta = \hatDta \cdot \M, \qquad \Eta = \hatEta \cdot \M.
\end{equation}

\section{Implementation in Matlab}

\label{s:matlab}

Given a regular enough function $f(t)$, in order to implement \eqref{e:DatEatfj}, we need to compute, on the one hand, the coefficients $\{\hat f_k\}$, for $0 \le k \le N$, and on the other hand, generate the matrices $\hatDta$ and $\hatEta$. We remark that we can safely operate with 64-bit precision to obtain $\{\hat f_k\}$, and that the final version of $\hatDta$ and $\hatEta$ can also be stored with 64-bit precision, as vpa is only necessary in the intermediate steps of the construction of $\hatDta$ and $\hatEta$. On the other hand, once that $\hatDta$ and $\hatEta$ have been obtained, it is straightforward to construct $\Dta$ and $\Eta$ from \eqref{e:tildeDattildeEat}, and then apply \eqref{e:tildeDatEatfj}. Note that $\M$ is a nonsingular well-conditioned matrix, and it can be safely generated with 64-bit precision, although using vpa to generate $\M\cdot\hatDta$ and $\M\cdot\hatEta$, and then convert the result into 64-bit precision yields higher accuracies.

\subsection{Obtention of the coefficients $\{f_k\}$}

The obtention of $\{f_k\}$ is a straightforward application of \eqref{e:gj}, \eqref{e:hatgk} and \eqref{e:gjhatgk}. Therefore, given the values $\{f(x_i)\}$, for $0\le j\le N$, we compute $g_j$, for $0\le j\le 2N-$, by means of \eqref{e:gjhatgk}, and finally, $\{\hat f_k\}$, for $0 \le N$, by means of \eqref{e:gjhatgk}:
\begin{equation*}
	\hat f_k = 
	\begin{cases}
		\hat g_0, & k = 0,
		\cr
		2\hat g_k, & 1 \le k \le N - 1,
		\cr
		\hat g_N, & k = N,
	\end{cases}
\end{equation*}
Note that the DFT in \eqref{e:gjhatgk} is performed by means of a fast Fourier transform (FFT) \cite{FFT}, that reduces the computational cost from $\mathcal O(N^2)$ to $\mathcal O(N\ln N)$. However, Matlab's implementation of the FFT multiplies the result by the number of input values, i.e., it actually calculates $\{2N\hat g_k\}$. Therefore, we have to divide the FFT result by $2N$:
$$
\hat g_k = \frac{1}{2N}\operatorname{FFT}(g_0, g_1, \ldots, g_{2N-1}) = \frac{1}{2N}\left[2N\sum_{j=0}^{2N-1}g_je^{-2ijk\pi/(2N)}\right].
$$
Bearing in mind the previous arguments, the Matlab code is as follows:
\begin{verbatim}
  g=[f;f(N:-1:2)];
  hatg=fft(g)/N;
  hatf=g(1:N+1);
  hatf([1 N+1])=hatf([1 N+1])/2;
  hatf(abs(hatf)<eps)=0;
\end{verbatim}
Observe that in the last line, we have applied the so-called Krasny's filter \cite{krasny}, i.e., we have set to zero those $f_k$ such that $|f_k| < \varepsilon$, where $\varepsilon = 2^{-52}$. Applying that filter is very common in FFT-based pseudospectral methods, and prevents the infinitesimally small rounding errors that appear when applying the FFT from getting amplified. Therefore, in general, we apply this filter systematically whenever it is possible.

\subsection{Computation of $\hatDta$ and $\hatEta$}

In order to generate efficiently $\hatDta$ and $\hatEta$, we observe that, thanks to the structure of \eqref{e:dj1k1} and \eqref{e:ej1k1}, they can be both expressed as a product of three matrices. More precisely, in the case of 
$\hatDta$:
\begin{equation}
\label{e:Dat}
\hatDta = \D_1\cdot\diag(\dd_2)\cdot\C(\lceil\alpha\rceil+1:N+1,1:N+1),
\end{equation}
where $\C(\lceil\alpha\rceil+1:N+1,1:N+1)\in\mathbb C^{(N+1-\lceil\alpha\rceil)\times(N+1)}$ is the matrix $\C$ without its first $\lceil\alpha\rceil$ columns,
\begin{equation*}
\D_1 = 
\begin{pmatrix}
t_0^{\lceil\alpha\rceil - \alpha} & t_0^{\lceil\alpha\rceil + 1 - \alpha} & \ldots & t_0^{N - \alpha}
	\cr
t_1^{\lceil\alpha\rceil - \alpha} & t_1^{\lceil\alpha\rceil + 1 - \alpha} & \ldots & t_1^{N - \alpha}
	\cr
\vdots & \vdots & \ddots & \vdots
	\cr
t_N^{\lceil\alpha\rceil - \alpha} & t_N^{\lceil\alpha\rceil + 1 - \alpha} & \ldots & t_N^{N - \alpha}	
\end{pmatrix}\in\mathbb C^{(N+1)\times(N+1-\lceil\alpha\rceil)},
\end{equation*}
and $\diag(\dd_2)\in\mathbb C^{(N+1-\lceil\alpha\rceil)\times(N+1-\lceil\alpha\rceil)}$ is the diagonal matrix whose diagonal entries are precisely the entries of the vector $\dd_2$:
\begin{equation}
\label{e:d2}
\dd_2 = \left(\frac{\Gamma(\lceil\alpha\rceil + 1)}{T^{\lceil\alpha\rceil}\Gamma(\lceil\alpha\rceil + 1 - \alpha)}, \frac{\Gamma(\lceil\alpha\rceil + 2)}{T^{\lceil\alpha\rceil+1}\Gamma(\lceil\alpha\rceil + 2 - \alpha)}, \ldots, \frac{\Gamma(N+1)}{T^N\Gamma(N + 1 - \alpha)}\right)\in\mathbb C^{N+1-\lceil\alpha\rceil}.
\end{equation}
Likewise, in the case of $\hatEta$,
\begin{equation}
\label{e:Eat}
\hatEta = \E_1\cdot\diag(\e_2)\cdot\C,
\end{equation}
where
\begin{equation*}
\E_1 = 
\begin{pmatrix}
t_0^{\alpha} & t_0^{1 + \alpha} & \ldots & t_0^{N + \alpha}
		\cr
t_1^{\alpha} & t_1^{1 + \alpha} & \ldots & t_1^{N + \alpha}
		\cr
\vdots & \vdots & \ddots & \vdots
		\cr
t_N^{\alpha} & t_N^{1 + \alpha} & \ldots & t_N^{N + \alpha}	
\end{pmatrix}\in\mathbb C^{(N+1)\times(N+1)},
\end{equation*}
and $\diag(\e_2)\in\mathbb C^{(N+1)\times(N+1)}$ is the diagonal matrix whose diagonal entries are precisely the entries of the vector $\e_2$:
\begin{equation}
\label{e:e2}
\e_2 = \left(\frac{\Gamma(1)}{T^{0}\Gamma(1 + \alpha)}, \frac{\Gamma(2)}{T^{1}\Gamma(2 + \alpha)}, \ldots, \frac{\Gamma(N+1)}{T^N\Gamma(N + 1 + \alpha)}\right)\in\mathbb C^{N+1}.
\end{equation}
With respect to the actual implementation in Matlab, we observe that the matrix $\C$ generated by \eqref{e:cij} is characterized by the fact that adjacent nonzero entries have opposite signs, and the absolute value of those entries can get extremely large values as $N$ grows. For instance, when $N = 100$, the exact values of $c_{71,101}$ and $c_{72,101}$ are
\begin{align*}
c_{71,101} & = 1666335331394399129847450283357894833563446572325067822868496043008453509120,
	\cr
c_{72,101} & = -1697794464111764171855358394790782869078728901979391848292743945528541184000,
\end{align*}
i.e., they both are of the order of $\mathcal O(10^{75})$ and have opposite signs. Therefore, working with standard 64-bit floating point precision (IEEE 754) causes most digits to be lost:
\begin{align*}
c_{71,101} & = 1.666335331394399\times10^{75},
	\cr
c_{72,101} & = -1.697794464111764\times10^{75},
\end{align*}
which has disastrous effects in the computation of \eqref{e:Dat} and \eqref{e:Eat}. In order to avoid this problem, we have found that the use of variable precision arithmetic (vpa), which allows working with arbitrary numbers of digits, allows computing \eqref{e:Dat} and \eqref{e:Eat} accurately. In Matlab, we have to fix the number of digits (by default, $32$) by means of the function \verb|digits|; then, the function \verb|vpa|, applied to a number, returns the number of digits required. For instance, if we type
\begin{verbatim}
  digits(1000)
  vpa(pi)
\end{verbatim}
the first $1000$ digits of $\pi$ are shown. Note that we could have also typed \verb|vpa(pi,1000)|, but if we have not set \verb|digits(1000)| previously, then, e.g., \verb|5*vpa(pi,1000)| returns only $32$ digits by default, so it is necessary to invoke \verb|digits| before using \verb|vpa| to maintain the desired accuracy along all the operations. On the other hand, sometimes it will be necessary to convert a given value into a symbolic variable, before applying \verb|vpa|. For instance, if we want $1000$ digits of $\pi^2$ and type \verb|vpa(pi^2,1000)|, we are not considering $\pi^2$, but its 64-bit floating point representation, so the result is incorrect. Therefore, we have to type commands like \verb|vpa(sym(pi)^2,1000)| or, even better, \verb|vpa(str2sym('pi^2'),1000)|, which turns the array of characters \verb|'pi^2'| into a symbolic variable. In general, the use of \verb|str2sym| is always a safe option. A final observation with respect to \verb|vpa| is that, if we want to apply it to an all-zero matrix, this is equivalent to applying \verb|sym|. Therefore, \verb|(vpa([0 0]) + 5).^1000 | and \verb|(sym([0 0]) + 5).^1000 | produce the same result, whereas \verb|(vpa(0) + 5)^100| and \verb|(sym(0) + 5)^1000| do not, because \verb|(vpa(0) + 5)^1000| returns the exact result, whereas \verb|(vpa(0) + 5)^1000| truncates the result to the number of digits set by the function \verb|digits|.

Bearing in mind the previous arguments, we generate first $\C$, which appears in both \eqref{e:Dat} and \eqref{e:Eat}, by implementing \eqref{e:cij} using \verb|vpa|. The resulting Matlab code is the following one:
\begin{verbatim}
  C=sym(zeros(N+1));
  C(1,1)=1;
  C(1:2,2)=[-1;2];
  for j=3:N+1
    C(1:j,j)=[(-1)^(j-1);4*C(1:j-1,j-1)-2*C(2:j,j-1)-C(2:j,j-2)];
  end
\end{verbatim}
Note that the only difference with respect to working with 64-bit floating point precision is the use of \verb|sym| in the first line. Indeed, if we type \verb|C=zeros(N+1);|, then, all the entries are generated with 64-bits, with the resulting loss of accuracy. On the other hand, we have used \verb|sym| in this piece of code, instead of \verb|vpa|, because, as we have said above, there is no difference in the result. If we want to impose that \verb|C| is really stored with \verb|vpa| and a specific number of digits, then, we must write, e.g., \verb|C(1,1)=vpa(1);| in the second line, but there is no significant variation in the elapsed time, so we have computed \verb|C| with \verb|sym|. However, even if it is possible to completely generate \eqref{e:Dat} and \eqref{e:Eat} by means of \verb|sym| without using \verb|vpa| at all, the computational cost is much more expensive, so we have only used \verb|sym| in the generation of \verb|C|.

The generation of $\D_1$ and $\E_1$ poses no problems. More precisely, after computing the time instants $\{t_j\}$ with \verb|vpa|,
\begin{verbatim}
  T_vpa=vpa(T);
  t=T_vpa*(1+cos(vpa(pi)*(0:N)/N)')/2;
\end{verbatim}
where \verb|T_vpa| contains the vpa version of the final time $T$; then, it is enough to type
\begin{verbatim}
  a_vpa=vpa(a);
  ceila=double(ceil(a_vpa));
  D1=t.^((ceila:N)-a_vpa);
  E1=t.^((0:N)+a_vpa);
\end{verbatim}
where \verb|a_vpa| contains  the vpa version of $\alpha$, and \verb|ceila| corresponds to $\lceil\alpha\rceil$.

In order to compute $\e_2$, we observe that
\begin{align*}
\frac{\Gamma(l + 1)}{\Gamma(l + 1 + \alpha)} = \frac{l}{l + \alpha}\cdot\frac{l - 1}{l - 1 + \alpha}\cdot\ldots\cdot\frac{1}{1 + \alpha}\cdot\frac{1}{\Gamma(1 + \alpha)}, \quad 0 \le l \le N,
\end{align*}
where we have used that $\Gamma(1 + z) = z\Gamma(z)$, for all $z\in\mathbb C$. Therefore, $\e_2$ in \eqref{e:e2} can be generated very efficiently as a comulative product:
\begin{verbatim}
  e2=cumprod([1/gamma(1+a_vpa),((1:N)./((1:N)+a_vpa))/T_vpa]).';
\end{verbatim}
observe that we store \verb|e2| as a column vector, because Matlab allows writing \verb|e2.*C|, which is equivalent to, but faster than \verb|diag(e2)*C|, so $\hatEta$ in \eqref{e:Eat} is created by typing
\begin{verbatim}
  hatEta=double(E1*(e2.*C));
\end{verbatim}
where we have converted the result from \verb|vpa| to \verb|double|, so \verb|hatEta| is stored in 64-bit precision. Observe that, even if we are asking $\alpha\not\in\mathbb N\cup\{0\}$ in the definition of $\hatDta$ and $\hatEta$, it is possible to consider also $\alpha\in\mathbb N\cup\{0\}$ to create \verb|hatEta|, without modifying the codes, but this is not true in the case of $\hatDta$. Indeed, when $\alpha\not\in\mathbb N\cup\{0\}$, we generate \verb|d2| in a similar way to \verb|e2|. More precisely, bearing in mind 
\begin{equation}
	\label{e:cumprodd2}
	\frac{\Gamma(l + 1)}{\Gamma(l + 1 - \alpha)} = \frac{l}{l - \alpha}\cdot\frac{l - 1}{l - 1 - \alpha}\cdot\ldots\cdot\frac{1}{1 - \alpha}\cdot\frac{1}{\Gamma(1 - \alpha)}, \quad 0 \le l \le N, 
\end{equation}
we write
\begin{verbatim}
  d2aux=cumprod([1/gamma(1-a_vpa),((1:N)./((1:N)-a_vpa))/T_vpa]).';
  d2=d2aux(ceila+1:N+1);
\end{verbatim}
i.e., we have ignored the first $\lceil\alpha\rceil$ entries of \verb|d2aux|. However, when $\alpha\in\mathbb N$, \eqref{e:cumprodd2} cannot be used, because $\Gamma(1 - \alpha)$ is not defined. Therefore, in that case, we write
\begin{equation*}
\frac{\Gamma(l + 1)}{\Gamma(l + 1 - \alpha)} = \frac{l}{l - \alpha}\cdot\frac{l - 1}{l - 1 - \alpha}\cdot\ldots\cdot\frac{2 + \alpha}{2}\cdot\frac{1 + \alpha}{1}\cdot\Gamma(1 + \alpha), \quad \alpha \le l \le N, 
\end{equation*}
and, hence, to generate $\dd_2$ in \eqref{e:d2}, we write
\begin{verbatim}
  d2=cumprod([prod(1:a_vpa)/T_vpa^a_vpa,((a_vpa+1):N)./(1:(N-a_vpa))/T_vpa]).';
\end{verbatim}
where \verb|prod(1:(1+a_vpa))| corresponds to $(1 + \alpha)\Gamma(1 + \alpha) = \Gamma(2 + \alpha) = (1 + \alpha)!$. Observe that, when $\alpha = 0$, \verb|d2| can be generated in both ways. Finally, in all cases, $\hatDta$ in \eqref{e:Dat} is created by typing
\begin{verbatim}
  hatDta=double(D1*(d2.*C(ceila+1:N+1,:)));
\end{verbatim}
where we have  removed the first $\lceil\alpha\rceil$ columns of $\C$, and have stored \verb|hatDta| in 64-bit precision.

\subsection{Computation of $\Dta$ and $\Eta$}

In order to obtain $\Dta$ and $\Eta$, we only need to compute $\M$ and apply \eqref{e:tildeDattildeEat}. We remark that it is numerically safe to construct $\M$ using 64-bit precision, and then multiply respectively the 64-bit versions of $\hatDta$ and $\hatEta$ by it. However, using vpa gives slightly more accurate results. The code is an immediate application of \eqref{e:mjk}:
\begin{verbatim}
  M=2*cos(vpa(pi)*(0:N)'*(0:N)/N);
  M(:,[1 N+1])=M(:,[1 N+1])/2;
  M([1 N+1],:)=M([1 N+1],:)/2;
  M=M/N;
\end{verbatim}
and the non-vpa version would follow by just replacing \verb|vpa(pi)| by \verb|pi|. Then, we multiply the vpa versions of \verb|hatDta| and \verb|hatEta| by it, and, afterwards, apply \verb|double| to all the matrices, i.e.,
\begin{verbatim}
  hatEta=E1*(e2.*C);
  Eta=double(hatEta*M);
  hatEta=double(hatEta);
\end{verbatim}
and
\begin{verbatim}
  hatDta=(D1*(d2.*C(ceila+1:N+1,:)));
  Dta=double(hatDta*M);
  hatDta=double(hatDta);
\end{verbatim}

\subsection{Matlab codes to generate $\hatDta$, $\hatEta$, $\Dta$ and $\Eta$}

In Listing~\ref{code:Dat}, we offer the code of the function \verb|CaputoMatrix.m|, that generates $\hatDta$ and $\Dta$, and in Listing~\ref{code:Eat}, the code of the function \verb|RiemannLiouvilleMatrix.m|, that generates $\hatEta$ and $\Eta$. Moreover, both functions return also the 64-bit discretization of $t\in[0,T]$. Even if, for the sake of clarity, we offer them separately, it is straightforward to combine them into a single function, which would imply some reduction of the execution time, because $\C$ and $\M$ would be then generated only once, instead of twice. In fact, if we increase the size of $\C$, new columns and rows are added, but the previously generated ones do not change. Therefore, it is perfectly possible to pregenerate a matrix $\C$ having a larger size than what is required, store it, and then introduce it as a parameter, a global variable, etc. In that case, we will have to type \verb|hatDta=D1*(d2.*C(ceila+1:N+1,1:N+1));| and \verb|hatEta=E1*(e2.*C(1:N+1,1:N+1));|, respectively.

Note also that the variables \verb|a_vpa| and \verb|T_vpa| store respectively a vpa version of the parameters \verb|a| (corresponding to $\alpha$) and \verb|T| (corresponding to the final time $T$). Therefore, when invoking \verb|CaputoMatrix.m| and \verb|RiemannLiouvilleMatrix.m|, it is possible to introduce \verb|a| and \verb|T| as 64-bit precision numbers, as vpa variables or even as symbolic variables; for instance, if $\alpha = \exp(-1)$, and $T = \pi^2$, then we must introduce \verb|a| and \verb|T| as \verb|str2sym('exp(-1)')| and \verb|str2sym('pi^2')|, respectively. On the other hand, in practice, it is not necessary to consider a vpa version of \verb|N| (which corresponds to $N$), because it contains a natural number that is correctly handled by Matlab.

\begin{lstlisting}[label=code:Dat, language=Matlab, basicstyle=\footnotesize, caption = {Matlab function \texttt{CaputoMatrix.m}, that generates $\hatDta$ and $\Dta$}]
% Create the Caputo differentiation matrices
% $\hat{\mathbf D}_t^\alpha$ and $\mathbf D_t^\alpha$.
% N, a, T and dig denote respectively
% $N$, $\alpha$, $\T$ and the number of digits for vpa.
function [hatDta,Dta,t]=CaputoMatrix(N,a,T,dig)
digits(dig); % Set the number of digits for vpa
a_vpa=vpa(a); % vpa version of $\alpha$
ceila=double(ceil(a_vpa)); % $\lceil\alpha\rceil$
T_vpa=vpa(T); % vpa version of $\T$
t=T_vpa*(1+cos(vpa(pi)*(0:N)/N)')/2; % $t\in[0,T]$
C=sym(zeros(N+1)); % $\mathbf C$
C(1,1)=1;
C(1:2,2)=[-1;2];
for j=3:N+1
    C(1:j,j)=[(-1)^(j-1);4*C(1:j-1,j-1)-2*C(2:j,j-1)-C(2:j,j-2)];
end
D1=t.^((ceila:N)-a_vpa); % $\mathbf D_1$
% $\mathbf d_2$
if double(a_vpa==ceil(a_vpa)) % Test whether $\alpha\in\mathbb Z$
    d2=cumprod([prod(1:a_vpa)/T_vpa^a_vpa,((a_vpa+1):N)./(1:(N-a_vpa))/T_vpa]).';
else
    d2aux=cumprod([1/gamma(1-a_vpa),((1:N)./((1:N)-a_vpa))/T_vpa]).'; 
    d2=d2aux(ceila+1:N+1);
end
hatDta=(D1*(d2.*C(ceila+1:N+1,:))); % $\hat{\mathbf D}_t^\alpha$
M=2*cos(vpa(pi)*(0:N)'*(0:N)/N);
M(:,[1 N+1])=M(:,[1 N+1])/2;
M([1 N+1],:)=M([1 N+1],:)/2;
M=M/N;
Dta=double(hatDta*M); % $\mathbf D_t^\alpha$
hatDta=double(hatDta);
t=double(t);
\end{lstlisting}

\begin{lstlisting}[label=code:Eat, language=Matlab, basicstyle=\footnotesize, caption = {Matlab function \texttt{RiemannLiouvilleMatrix.m}, that generates $\hatEta$ and $\Eta$}]
% Create the Riemann-Liouville integration matrices
% $\hat{\mathbf E}_t^\alpha$ and $\mathbf E_t^\alpha$.
% N, a, T and dig denote respectively
% $N$, $\alpha$, $\T$ and the number of digits for vpa.
function [hatEta,Eta,t]=RiemannLiouvilleMatrix(N,a,T,dig)
digits(dig); % Set the number of digits for vpa
a_vpa=vpa(a); % vpa version of $\alpha$
T_vpa=vpa(T); % vpa version of $\T$
t=T_vpa*(1+cos(vpa(pi)*(0:N)/N)')/2; % $t\in[0,T]$
C=sym(zeros(N+1)); % $\mathbf C$
C(1,1)=1;
C(1:2,2)=[-1;2];
for j=3:N+1
    C(1:j,j)=[(-1)^(j-1);4*C(1:j-1,j-1)-2*C(2:j,j-1)-C(2:j,j-2)];
end
E1=t.^((0:N)+a_vpa); % $\mathbf E_2$
e2=cumprod([1/gamma(1+a_vpa),((1:N)./((1:N)+a_vpa))/T_vpa]).'; % $\mathbf e_2$
hatEta=E1*(e2.*C); % $\hat{\mathbf E}_t^\alpha$
M=2*cos(vpa(pi)*(0:N)'*(0:N)/N);
M(:,[1 N+1])=M(:,[1 N+1])/2;
M([1 N+1],:)=M([1 N+1],:)/2;
M=M/N;
Eta=double(hatEta*M); % $\mathbf E_t^\alpha$
hatEta=double(hatEta);
t=double(t);
\end{lstlisting}

\subsection{Estimation of the number of digits}

In general, except for very small values of $N$, it is not possible to generate $\hatDta$, $\Dta$, $\hatEta$ and $\Eta$ correctly, a fact that makes compulsory the use of vpa. It is easy to check this experimentally, by just removing the variable \verb|dig| and all the appearences of \verb|vpa| in \verb|CaputoMatrix.m| (see Listing~\ref{code:Dat}) and \verb|RiemannLiouvilleMatrix.m| (see Listing~\ref{code:Eat}). Then, if we take $N = 40$, $\alpha = 0.37$ and $T = 1.2$, the thus generated matrices $\hatDta$, $\Dta$, $\hatEta$ and $\Eta$ are nonsensical and contain entries whose moduli are as large as $\mx{\hatDta} = 2.4053\times10^{14}$, $\mx{\Dta} = 8.2056\times10^{12}$, $\mx{\hatEta} = 2.6064\times10^{13}$ and $\mx{\Eta} = 8.8991\times10^{11}$, respectively, where $\mx{\cdot}$ denotes the maximum norm (recall that, given a matrix $\A = [a_{ij}]$, $\mx{\A} = \max_{i,j}|a_{ij}|$ returns the maximum of  the moduli of all the entries of $\A$).

This behavior is also appreciated if we work with vpa and a too small number of digits \verb|dig|. For instance, if $N = 100$, $\alpha = 0.37$ and $T = 1.2$, then, $\mx{\hatDta(\dig=30)} = 1.7477\times10^{38}$, $\mx{\Dta(\dig=30)} = 3.8319\times10^{36}$, $\mx{\hatEta(\dig=30)} = 3.3523\times10^{37}$ and $\mx{\Eta(\dig=30)} = 3.7454\times10^{35}$, which makes evident that $30$ digits are not enough when $N = 100$. However, if we increase the number of digits to $100$, then $\mx{\hatDta(\dig=100)} = 46.0508$, $\mx{\Dta(\dig=100)} = 26.2840$, $\mx{\hatEta(\dig=100)} = 1.2029$ and $\mx{\Eta(\dig=100)} = 0.19984$. Moreover, once that the necessary number of digits \verb|dig| has been chosen, there is no difference in the final 64-bit precision version of the matrices, if we further increase \verb|dig|. For instance, after choosing 1000 digits, the results are identical to those with 100 digits. Therefore, given an adequate value of \verb|dig|, we can simply set $\hatDta = \hatDta(\dig)$, $\Dta = \Dta(\dig)$, $\hatEta = \hatEta(\dig)$ and $\Eta = \Eta(\dig)$.

\begin{figure}[!htbp]
	\centering
	\includegraphics[width=0.5\textwidth, clip=true]{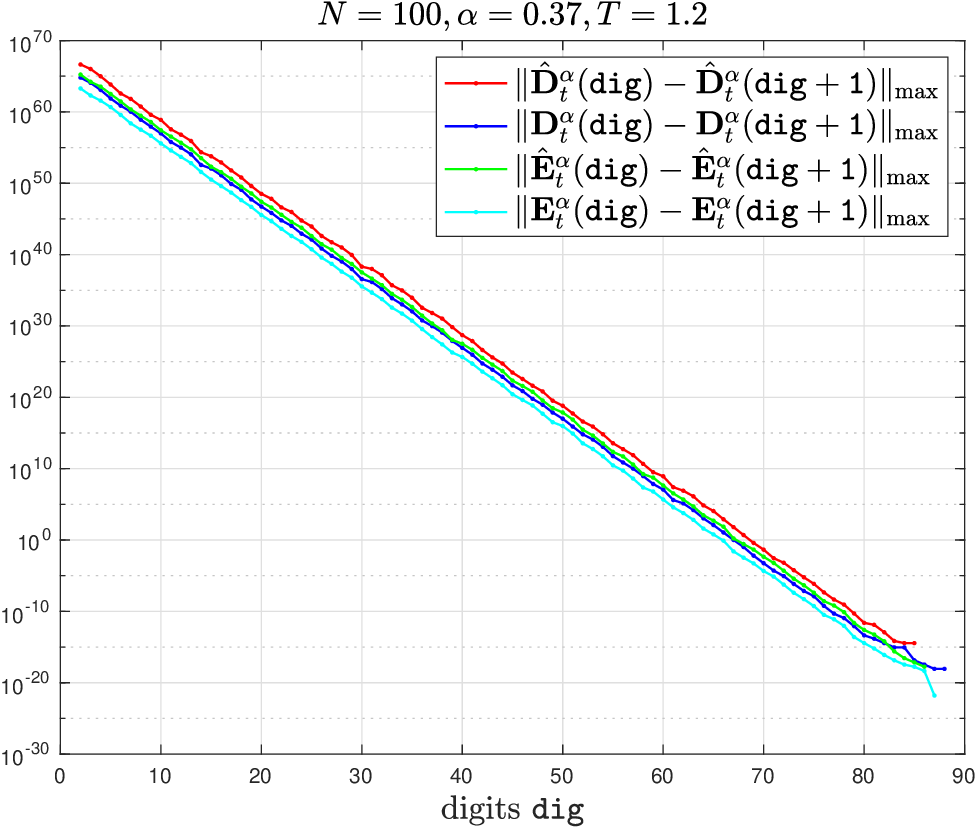}
	\caption{$\mx{\hatDta(\dig) - \hatDta(\dig+1)}$ (red), $\mx{\Dta(\dig) - \Dta(\dig+1)}$ (blue), $\mx{\hatEta(\dig) - \hatEta(\dig+1)}$ (green) and $\mx{\Eta(\dig) - \Eta(\dig+1)}$ (cyan), for $\mathtt{dig} = 2, 3, 4, \ldots$. In all cases, we have taken $N = 100$, $\alpha = 0.37$ and $T = 1.2$.}
	\label{f:mindigDE}
\end{figure}

In Figure~\ref{f:mindigDE}, in order to clarify the effect of \verb|dig| on the generation of the matrices, we have plotted in semilogarithmic scale $\mx{\hatDta(\dig) - \hatDta(\dig+1)}$, $\mx{\Dta(\dig) - \Dta(\dig+1)}$, $\mx{\hatEta(\dig) - \hatEta(\dig+1)}$ and $\mx{\Eta(\dig) - \Eta(\dig+1)}$ with respect to $\mathtt{dig} = 2, 3, 4, \ldots$, until these quantities are strictly equal to zero (although it is enough to set that they are smaller than a given infinitesimal number $\varepsilon$, e.g., $\varepsilon = 2^{-52}$, which is given in Matlab by \verb|eps|). As can be seen, they decay exponentially with respect to \verb|dig|.

\begin{figure}[!htbp]
	\centering
	\includegraphics[width=0.5\textwidth, clip=true]{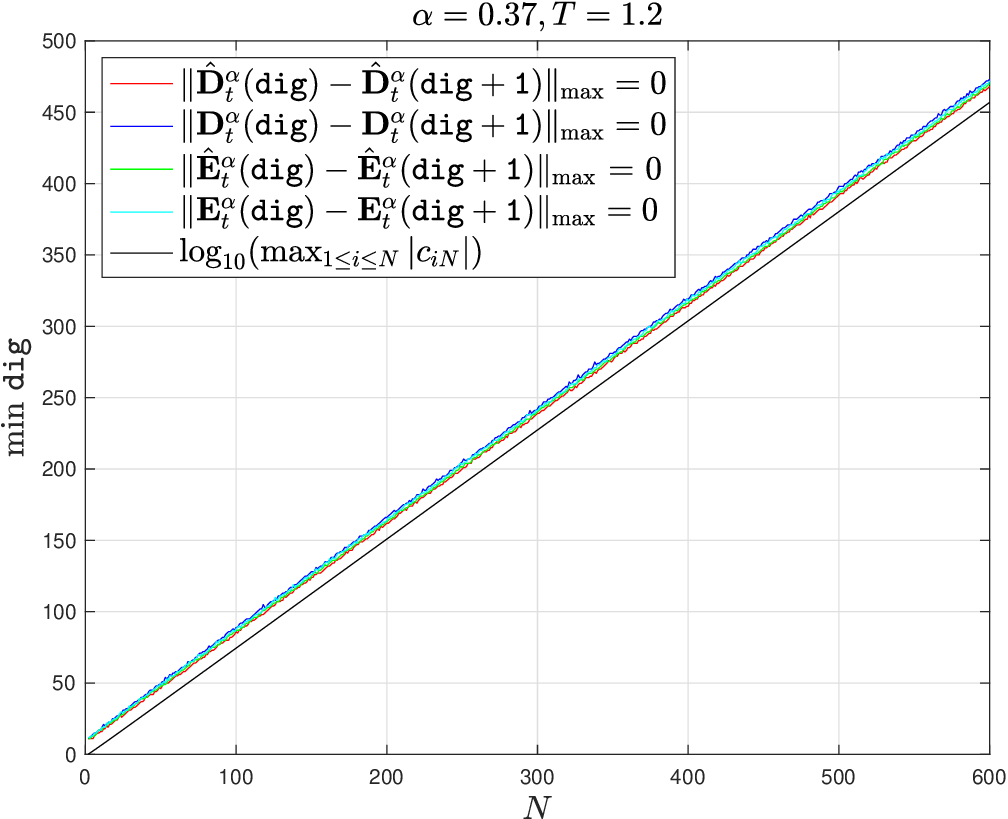}\includegraphics[width=0.5\textwidth, clip=true]{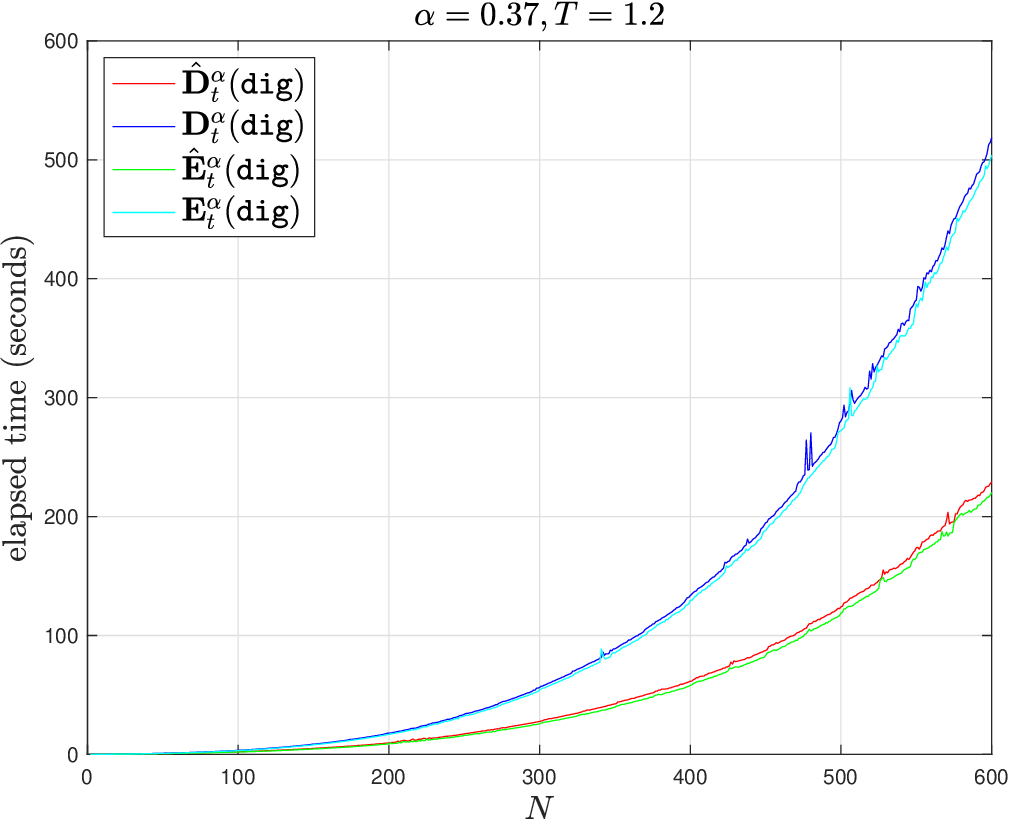}
	\caption{Left: minimum values of \texttt{dig} for which both $\mx{\hatDta(\dig) - \hatDta(\dig+1)} = 0$ and $\mx{\Dta(\dig) - \Dta(\dig+1)} = 0$ (red), and minumum values of \texttt{dig} for which both $\mx{\hatEta(\dig) - \hatEta(\dig+1)} = 0$ and $\mx{\Eta(\dig) - \Eta(\dig+1)} = 0$ (blue), together with $\max_{1\le i\le N}|c_{iN}|$, for $2\le N \le 500$. Right: elapsed time of the experiments on the left-hand side. In all cases, we have taken $\alpha = 0.37$ and $T = 1.2$.}
	\label{f:mindigNDE}
\end{figure}

To better understand how to choose \verb|dig| in terms of $N$, we show, on the left-hand side of Figure~\ref{f:mindigNDE}, the minumum values of \verb|dig| for which $\mx{\hatDta(\dig) - \hatDta(\dig+1)} = 0$ (red), $\mx{\Dta(\dig) - \Dta(\dig+1)} = 0$ (blue), $\mx{\hatEta(\dig) - \hatEta(\dig+1)} = 0$ (green) and $\mx{\Eta(\dig) - \Eta(\dig+1)} = 0$ (cyan), for $2\le N\le 600$;  these conditions are very exigent and, as mentioned above, can be relaxed, by setting, e.g., $< \varepsilon=2^{-52}$ instead of $= 0$. The numerical experiments reveal that the number of necessary digits grows up approximately linearly with $N$, and a linear regression analysis yields $\min\dig \approx 0.7663N + 8.449$ for $\hatDta$, $\min\dig \approx 0.7704N + 11.57$ for $\Dta$, $\min\dig \approx 0.7671N + 9.59$ for $\hatEta$ and $\min\dig \approx 0.7684N + 10.81$ for $\Eta$.  For the sake of comparison, we have also plotted on the left-hand side of Figure~\ref{f:mindigNDE} the decimal logarithm of the maximum of the absolute values of the $N$th column of $\C$ in terms of $N$, i.e., $\log_{10}(\max_{1\le i\le N}|c_{iN}|) \approx 0.7645N - 1.946$, i.e., the fifth regression lines are approximately parallel, which shows that the number of necessary digits is correlated with the number of digits necessary to store the corresponding column of $\C$ in an exact way. 

On the other hand, we have plotted on the right-hand side of Figure~\ref{f:mindigNDE} the times needed to generate $\hatDta(\dig)$ (red), $\Dta(\dig)$ (blue), $\hatEta(\dig)$ (green) and $\Eta(\dig)$ (cyan), for the corresponding values of $\dig$ obtained on the left-hand side of Figure~\ref{f:mindigNDE}; note that, in the case of $\hatDta$ and $\hatEta$, we have executed \verb|CaputoMatrix.m| and \verb|RiemannLiouvilleMatrix.m|, after removing the code corresponding to the generation of $\M$, and of $\Dta$ and $\Eta$, which diminishes the execution time, because the computaton of $\Dta$ and $\Eta$ requires $\hatDta$ and $\hatEta$, respectively, but not the other way around. In all cases, the number of seconds seems to grow as $\mathcal O(N^3)$, and a least-square regression analysis yields $\mathrm{time}\approx 1.568\times10^{-6}N^3 - 0.000454N^2 + 0.09306N - 2.875$ for $\hatDta$, $\mathrm{time}\approx 3.59\times10^{-6}N^3 - 0.001004N^2 + 0.1744N - 5.329$ for $\Dta$, $\mathrm{time}\approx 1.56\times10^{-6}N^3 - 0.0004728N^2 + 0.0924N - 2.858$ for $\hatEta$, and $\mathrm{time}\approx 3.585\times10^{-6}N^3 - 0.001041N^2 + 0.1789N - 5.557$ for $\Eta$.

Finally, let us mention that the equivalent results to those in Figure~\ref{f:mindigNDE} are very similar for other values of $\alpha$ and $T$, but, in case of doubt, it is enough to generate the required matrices for \verb|dig| and $\texttt{dig}+1$, and check whether they are identical.

\section{Numerical tests}

\label{s:numerical}

We have used the functions \verb|CaputoMatrix.m| and  \verb|RiemannLiouvilleMatrix.m| written in Section~\ref{s:matlab}, to approximate $\cpt e^{imt}$ and $\rl e^{imt}$ respectively, whose exact expressions can be obtained, e.g, by typing in Mathematica \cite{mathematica} the following:
\begin{verbatim}
	f[t_]=Exp[I m t];
	Integrate[D[f[tau],{tau,n}]/(t-tau)^(1-n+a),{tau,0,t}]/Gamma[n-a]
	Integrate[f[tau]/(t-tau)^(1-a),{tau,0,t}]/Gamma[a]
\end{verbatim}
which, after some minor simplification, yields
\begin{align}
\label{e:cpteimt}
	\cpt e^{imt} & = (im)^\alpha e^{imt}\left(1 - \frac{\Gamma(\lceil\alpha\rceil - \alpha, imt)}{\Gamma(\lceil\alpha\rceil - \alpha)}\right), \quad \alpha\not=\lceil\alpha\rceil,
	\\
\label{e:rleimt}
	\rl e^{imt} & = (im)^{-\alpha} e^{imt}\left(1 - \frac{\Gamma(\alpha, imt)}{\Gamma(\alpha)}\right), \quad \alpha > 0,
\end{align}
where $\Gamma(\cdot, \cdot)$, whose corresponding Matlab command is \verb|igamma|, denotes the upper incomplete gamma function:
$$
\Gamma(s, x) = \int_{x}^{\infty}t^{s-1}e^{-t}dt,
$$
and hence, $\Gamma(s) = \Gamma(s, 0)$.

\subsection{A full numerical test in Matlab}

In Listing~\ref{code:test}, we offer a Matlab program that approximates numerically $\cpt e^{imt}$ and $\rl e^{imt}$. In the case of $\cpt e^{imt}$, we compare the results with \eqref{e:cpteimt}, when $\alpha\not=\lceil\alpha\rceil$, and with $(d^\alpha/dt^\alpha)e^{imt} = (im)^\alpha e^{imt}$, when $\alpha\not=\lceil\alpha\rceil$, in particular, $\alpha = 0$ returns $e^{imt}$ itself. On the other hand, in the case of $\rl e^{imt}$, we compare the results with \eqref{e:rleimt}, when $\alpha > 0$, and with $e^{imt}$ itself, when $\alpha = 0$. For the sake of comparison, we have also generated $\M$, as in \eqref{e:mjk}, using 64-bit precision, and compared $\hatDta\cdot\M$ and $\hatEta\cdot\M$ to $\Dta$ and $\Eta$, respectively.

The experiment for $\alpha = 1.3\in(1,2)$, $T = 1.2$, $m = 2$, $N = 100$ and $\dig=100$ digits takes 8.16 seconds to run, yielding the following errors in discrete $\ell^\infty$ norm:
\begin{itemize}
	\item $\displaystyle\max_{0\le j\le N}\big|[\M\cdot\f]_j - [\hatf]_j\big| = 3.5376\times10^{-15}$.
	
	\item $\displaystyle\max_{0\le j\le N}\big|[\hatDta\cdot\hatf]_j - \cpt f(t_j)\big| =
	\left\{\begin{aligned}
	& 2.8893\times10^{-11} & & \mbox{(without Krasny's filter)},
		\cr
	& 6.8315\times10^{-14} & & \mbox{(with Krasny's filter)}.
	\end{aligned}\right.$
	
	\item $\displaystyle\max_{0\le j\le N}\big|[\Dta\cdot\f]_j - \cpt f(t_j)\big| = 3.7006\times10^{-11}$.

	\item $\displaystyle\max_{0\le j\le N}\big|[\hatDta\cdot\M\cdot\f]_j - \cpt f(t_j)\big| = 1.0433\times10^{-9}$.
	
	\item $\displaystyle\max_{0\le j\le N}\big|[\hatEta\cdot\hatf]_j - \rl f(t_j)\big| = 3.5108\times10^{-16}$ (error with and without Krasny's filter).

\item $\displaystyle\max_{0\le j\le N}\big|[\Eta\cdot\f]_j - \rl f(t_j)\big| = 4.5776\times10^{-16}$.

\item $\displaystyle\max_{0\le j\le N}\big|[\hatEta\cdot\M\cdot\f]_j - \rl f(t_j)\big| = 4.7429\times10^{-16}$.

\end{itemize}

As can be seen, the errors in the approximation of $\rl f(t)$ are always of the order of $10^{-16}$, and Krasny's filter is unnecessary. However, when approximating $\cpt f(t)$, $\hatDta\cdot\hatf$ gives a better approximation than $\Dta\cdot\f$, because of Krasny's filter. On the other hand, even if $\M\cdot\f$ calculates $\hatf$ with great accuracy, the worst results are obtained with $\hatDta\cdot\M\cdot\f$. Finally, let us mention that the choice of \verb|dig| has a limited impact the time of execution. For instance, if we take \verb|dig| four times as large, i.e., $\dig=400$, then the elapsed time is of 9.47 seconds.

\begin{lstlisting}[label=code:test, language=Matlab, basicstyle=\footnotesize, caption = {Matlab program \texttt{testmatrices.m}, that approximates numerically $\cpt e^{it}$ and $\rl e^{it}$}]
clear
tic
format short
a=0.37; % $\alpha$
N=100; % $N$
dig=100; % Number of digits
T=1.2; % $T$
% $\hat{\mathbf D}^\alpha_t$, $\mathbf D^\alpha_t$ and $t$
[hatDta,Dta,t]=CaputoMatrix(N,a,T,dig);
% $\hat{\mathbf E}^\alpha_t$ and $\mathbf E^\alpha_t$
[hatEta,Eta]=RiemannLiouvilleMatrix(N,a,T,dig);
m=2;
f=exp(1i*m*t); % $f(t)=e^{imt}$
% Exact values of $D^\alpha_t f(t_j)$
if ceil(a)==a % if $\alpha$ is an integer
    Dtaf=(1i*m)^a*exp(1i*m*t);
else
    Dtaf=(1i*m)^a*exp(1i*m*t).*(1-igamma(ceil(a)-a,1i*m*t)/gamma(ceil(a)-a));
end
% Exact values of $I^\alpha_t f(t_j)$
if a == 0 % if $\alpha=0$
    Itaf=exp(1i*m*t);
else
    Itaf=(1i*m)^(-a)*exp(1i*m*t).*(1-igamma(a,1i*m*t)/gamma(a));
end
g=[f;f(N:-1:2)]; % $g_j$
hatg=fft(g)/N; % $\hat g_k$
hatf=hatg(1:N+1); % $\hat f_k$
hatf([1 N+1])=hatf([1 N+1])/2;
M=2*cos(pi*(0:N)'*(0:N)/N); % Generate $M$ using 64-bit precision
M(:,[1 N+1])=M(:,[1 N+1])/2;
M([1 N+1],:)=M([1 N+1],:)/2;
M=M/N;
% Errors in $L^\infty$ norm
norm(hatf-M*f,inf) % error in the approximation of $\hat f$ by $M$
norm(hatDta*hatf-Dtaf,inf) % error without Krasny's filter
norm(hatEta*hatf-Itaf,inf) % error without Krasny's filter
hatf(abs(hatf)<eps)=0; % Krasny's filter
norm(hatDta*hatf-Dtaf,inf) % error with Krasny's filter
norm(hatEta*hatf-Itaf,inf) % error with Krasny's filter
norm(Dta*f-Dtaf,inf)
norm(hatDta*M*f-Dtaf,inf)
norm(Eta*f-Itaf,inf)
norm(hatEta*M*f-Itaf,inf)
toc % Elapsed time
\end{lstlisting}

\subsection{Numerical approximation of $\cpt e^{imt}$ and $\rl e^{imt}$ for large natural values of $m$ and different values of $\alpha$ and $N$}

One important application of the matrices that we have constructed is the numerical approximation of the Caputo derivative and Riemann-Liouville integral of highly oscillatory integrals. In this regard, we have considered $\cpt e^{imt}$ and $\rl e^{imt}$, taking $m = 110$, for which a large amount of nodes $N$ is needed. The numerical experiments reveal that, in general, the approximation of $\rl e^{imt}$ poses no problems and spectrally accurate results are always obtained, provided that $N$ is large enough. On the other hand, due to the term $(im)^\alpha$ in \eqref{e:cpteimt}, $\max_{t\in[0,T]}|\cpt e^{imt}|$ grows rapidly with $\alpha$; therefore, we have used the maximum relative error, which enables comparing the results for different values of $\alpha$ in a more consistent way. More precisely, given $f(t) = e^{imt}$, we define
\begin{equation}
	\label{e:errhatDta}
	\errhatDta = 
\left\{
\begin{aligned}
& \max_{0\le j\le N-1}\left|\frac{[\hatDta\cdot\hatf]_j - \cpt f(t_j)}{\cpt f(t_j)}\right|, & & \alpha\not=\lceil\alpha\rceil,
	\cr
& \max_{0\le j\le N}\left|\frac{[\hatDta\cdot\hatf]_j - (im)^\alpha e^{imt}}{(im)^\alpha e^{imt}}\right|, & & \alpha=\lceil\alpha\rceil,
\end{aligned}
\right.	
\end{equation}
and
\begin{equation}
	\label{e:errDta}
	\errDta = 
	\left\{
	\begin{aligned}
		& \max_{0\le j\le N-1}\left|\frac{[\Dta\cdot\hatf]_j - \cpt f(t_j)}{\cpt f(t_j)}\right|, & & \alpha\not=\lceil\alpha\rceil,
		\cr
		& \max_{0\le j\le N}\left|\frac{[\Dta\cdot\hatf]_j - (im)^\alpha e^{imt}}{(im)^\alpha e^{imt}}\right|, & & \alpha=\lceil\alpha\rceil,
	\end{aligned}
	\right.	
\end{equation}
where $\cpt f(t_j)$ is given by \eqref{e:cpteimt}. Note that, when $\alpha\not=\lceil\alpha\rceil$, we have omitted the value $j = N$ corresponding to $t_N = 0$, because $\cpt f(0) = 0$ in that case.

\begin{figure}[!htbp]
	\centering
	\includegraphics[width=0.5\textwidth, clip=true]{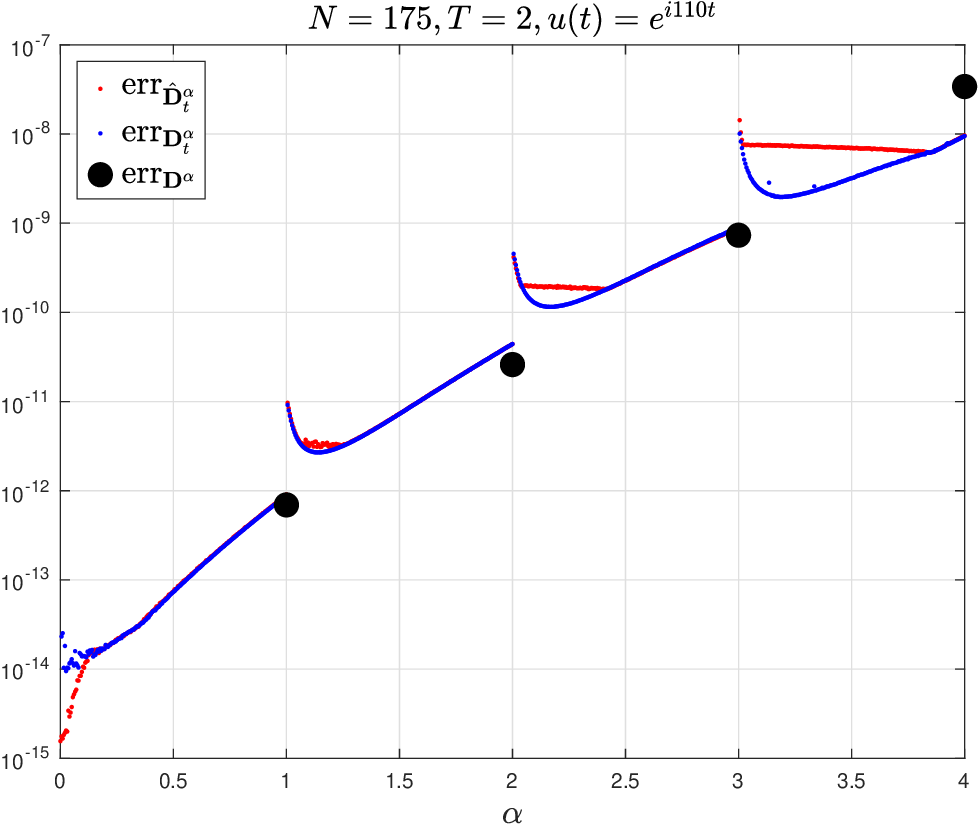}\includegraphics[width=0.5\textwidth, clip=true]{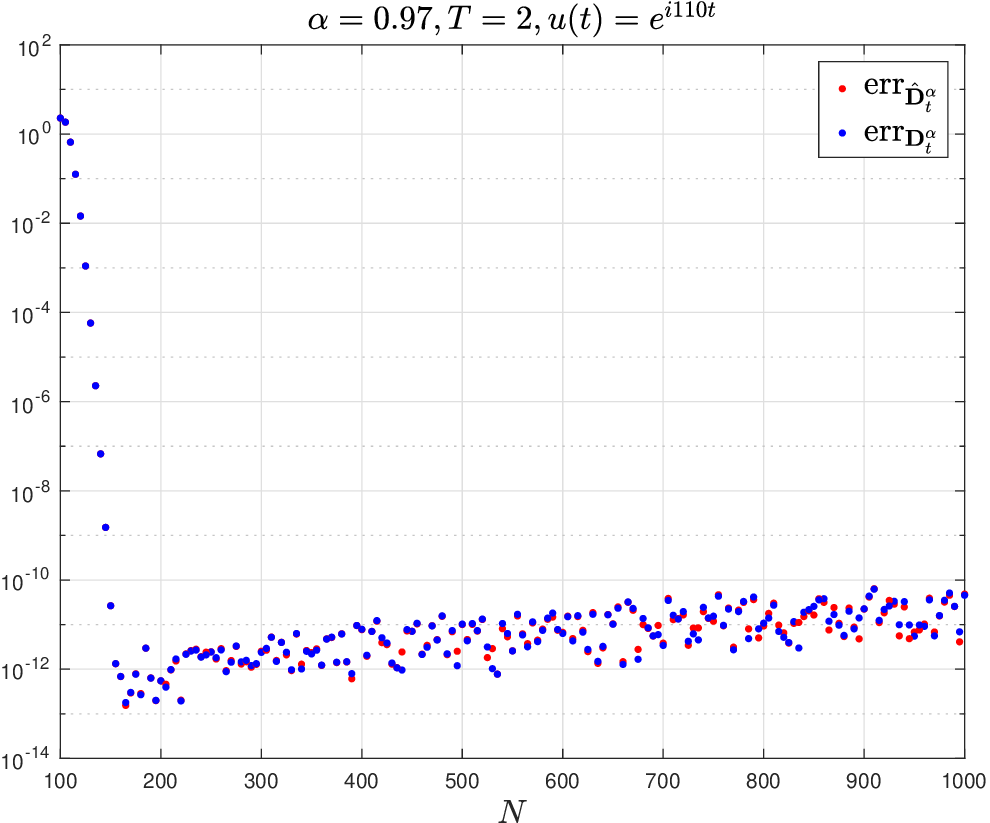}
	\caption{Left: Errors $\errhatDta$ (red) and $\errDta$ (blue) as defined in \eqref{e:errhatDta} and \eqref{e:errDta}, respectively, for $f(t) = e^{i110t}$, i.e., $m = 110$, and $\alpha\{0, 0.005, \ldots, 4\}$. For the sake of comparison, we also plot $\err_{\D^\alpha}$, as defined in \eqref{e:errhatDa}, for $\alpha\in\{1, 2, 3, 4\}$ (thick black point). Right: Errors $\errhatDta$ (red) and $\errDta$ (blue), for $f(t) = e^{i110t}$, and $N\in\{100, 105, \ldots, 1000\}$.}
	\label{f:errNN}
\end{figure}

Since, when $\alpha\in\mathbb N\cup\{0\}$, our numerical approximation of $\cpt$ approximates the standard integer-order derivative of order $\alpha$, it follows that, when, e.g., $\Dta$ should offer similar results to those given by the standard Chebyshev-differentiation matrix $\D$ generated by \verb|cheb.m| \cite{trefethen}. Therefore, in order to do a fair comparison between $\Dta$ and $\D$, we have taken $T = 2$, so $t_0 - t_N = 2$ and $\D$ does not need to be scaled. On the left-hand side of Figure~\ref{f:errNN}, we have plotted $\errhatDta$ (red) and $\errDta$ (blue), for $N = 175$, $T = 2$ and $\alpha\in\{0, 0.005, 0.01, \ldots, 4\}$, taking $\mathtt{dig} = N$, and, additionally, we have plotted (thick black point), for $\alpha\in\{1,2,3,4\}$,
\begin{equation}
\label{e:errhatDa}
\err_{\D^\alpha} = \max_{0\le j\le N}\left|\frac{[\D^\alpha\cdot\hatf]_j - (im)^\alpha e^{imt}}{(im)^\alpha e^{imt}}\right|,
\end{equation}
where $\D^\alpha$ denotes the power $\alpha$ of $\D$. The results show that the errors corresponding to $\hatDta$ and $\Dta$ are very similar, and the same is valid for the errors corresponding to $\D$ and natural values of $\alpha$. Indeed, as with $\D^\alpha$, the results get worse as $\alpha$ increases, but they are remarkably good for $\alpha\in[0, 2]$ and, especially, for $\alpha\in[0, 1]$. Note however that the graphic of the errors exhibits jumps at the natural values of $\alpha$, because, the operator $\cpt$ is not continuous at those values of $\alpha$. On the other hand, in what regards the numerical approximation of $\rl e^{imt}$, the errors in discrete $\ell^\infty$ norm are always of the order of $\mathcal O(10^{-15})$, for all the values of $\alpha$ considered.

In order to see the effect of $N$ on the errors, we have plotted, on the right-hand side of Figure~\ref{f:errNN}, $\errhatDta$ (red) and $\errDta$ (blue), for $\alpha = 0.97$, i.e., a value close to $1$, and $N \in\{100, 105, \ldots, 1000\}$, taking $\dig = N$, which is a more than enough amount of digits. The results are almost identical for $\errhatDta$ and $\errDta$, and reveal that, even if these quantities show a weak tendency to grow with $N$, they always stay lower than $10^{-10}$. Therefore, we can say that the generation of $\hatDta$ and $\Dta$ remains stable, even when large values of $N$ are considered.  In what regards the numerical approximation of $\rl e^{imt}$, the errors in discrete $\ell^\infty$ norm are again of the order of $\mathcal O(10^{-15})$, in the numerical experiments with $N\ge 155$.

\section{Numerical simulation of Caputo-type advection-diffusion equations}

\label{s:advdif}

The developed matrices allow us, as in \cite{delahozmuniain2024}, to approximate numerically Caputo-type advection-diffusion equations of the form
\begin{equation}
	\label{e:cptulinear}
	\begin{cases}
		\cpt u(t,x) = a_1(x)u_{xx}(t,x) + a_2(x)u_x(t,x) + a_3u(t,x) + a_4(t,x), \quad t\in[0, T],
		\cr
		u(x, 0) = u_0(x),
	\end{cases}
\end{equation}
where $a_1(x)$, $a_2(x)$, $a_3(x)$, $a_4(t,x)$ are at least continuous, and $x\in\mathbb R$ or $x\in[x_a, x_b]$, although other types of spatial domains like semi-infinite ones can be considered, as well as PDEs containing higher order derivatives with respect to $x$. As said in the introduction, the idea is to transform \eqref{e:cptulinear} into a Sylvester equation of the form $\A\cdot\X+\X\cdot\B=\C$, then solve it. Moreover, by using the ideas in, e.g, \cite{cuestadelahoz2024}, we can consider generalizations of \eqref{e:cptulinear} to $d$ spatial dimensions, with $d > 1$, and transform them to the so-called Sylvester tensor equations. Therefore, let us consider first an example in one spatial dimension, then another example in multiple spatial dimensions.

\subsection{An example in one spatial dimension}

\label{s:oneD}

In one spatial dimension, the procedure is very similar to that explained in \cite{delahozmuniain2024}, so we omit many details, which can be consulted in that reference, especially how to incorporate the boundary conditions. Note however that, unlike in \cite{delahozmuniain2024}, the nodes $t_j\in[0, T]$, as defined in \eqref{e:tj}, are given in decreasing order, i.e., $t_0 = T > t_1 > \ldots > t_{N} = 0$, which implies a small modification of the codes in \cite{delahozmuniain2024}. To illustrate this, we consider the following example:
\begin{equation}
	\label{e:PDE1}
	\left\{
	\begin{aligned}
		& \cpt u(t,x) = u_{xx}(t,x) + 2xu_x(t, x) + 2u(t,x)
		\cr
		& \qquad\qquad\qquad + (im)^\alpha\left(1 - \frac{\Gamma(\lceil\alpha\rceil - \alpha, imt)}{\Gamma(\lceil\alpha\rceil - \alpha)}\right)e^{imt - x^2}, & & (t,x)\in[0, 2]\times\mathbb R,
		\cr
		& u(x, 0) = u_0(x) = e^{-x^2},
	\end{aligned}
	\right.
\end{equation}
where $\alpha\in(0,1)$, $m\in\N$, and whose solution is $u(t,x) = e^{imt - x^2}$. In this case, we discretize the operator $\cpt$ by means of $\Dta$, and the first and second spatial derivatives by means of $\Dx$ and $\Dx^2$, where $\Dx\in\R^{N_x\times N_x}$ is the Hermite differentiation matrix generated by the Matlab function \verb|herdif| (see \cite{WeidemanReddy2000}), which is based on the Hermite functions
$$
\psi_n(x) = \frac{e^{-x^2/2}}{\pi^{1/4}\sqrt{2^nn!}}H_n(x),
$$
where $H_n(x)$ are the Hermite polynomials:
$$
H_n(x) = \frac{(-1)^n}{2^n}e^{x^2}\frac{d^n}{dx^n}e^{-x^2},
$$
and the spatial nodes $x_j$, with $0 \le j \le N_x-1$, are the roots of the polynomial $H_n(x)$ with the largest index, multiplied by a scale factor $b$, which is precisely the third parameter of \verb|herdif|. Therefore, we want to obtain a matrix $\U = [u_{jk}]\in\CC^{(N_t+1)\times N_x}$, such that $u_{jk} \approx u(t_j, u_k)$; note that we write  $N_t$ instead of $N$, to avoid any confusion with $N_x$.

Following the ideas in \cite{delahozmuniain2024}, let us denote
\begin{equation}
\label{e:EF}
\mathbf E =
\left(
\begin{array}{ccc}
	1 & & 0
	\cr
	& \ddots &
	\cr
	0 & & 1
	\cr
\hline
	0 & \ldots & 0
\end{array}
\right),
\qquad
\mathbf F =
\left(
\begin{array}{ccc}
	0 & \ldots & 0
	\cr
	\vdots & \ddots & \vdots
	\cr
	0 & \ldots & 0
	\cr
	\hline
		u_0(x_0) & \ldots & u_0(x_{N_x})
\end{array}
\right);
\end{equation}
note that, due to the fact that the nodes $t_j$ are defined in decreasing order, we store an all-zero vector in the last row of $\E$, and the initial condition $u_0(x)$ is stored in the last row of $\F$, whereas in \cite{delahozmuniain2024}, they were stored in the first rows of $\E$ and $\F$, respectively. Then, we decompose $\U$ as
\begin{equation}
	\label{e:Uboundaryt}
	\U = \E\cdot\Uin+ \F,
\end{equation}
where $\Uin\in\CC^{N_t\times N_x}$ denotes the matrix of $\U$ without its last row, but keeping all its columns (in \cite{delahozmuniain2024}, we do not keep the first row of $\U$ to define $\Uin$):
\begin{equation*}
	\Uin =
	\begin{pmatrix}
		u_{11} & \ldots & u_{1,N_x}
		\cr
		\vdots & \ddots & \vdots
		\cr
		u_{N_t,1} & \ldots & u_{N_t,N_x}
	\end{pmatrix}.
\end{equation*}
Moreover, denoting $\x = (x_0, x_1, \ldots, x_{N_x-1})$, and $\diag(\x)$ the diagonal matrix whose diagonal entries are precisely those of $\x$, we define
\begin{equation*}
\G = \Dx^2 + 2\diag(\x)\Dx + 2\I\in\R^{N_x\times N_x},
\end{equation*}
and $\HH = [h_{jk}]\in C^{(N_t+1)\times N_x}$, such that
$$
h_{jk} = (im)^\alpha\left(1 - \frac{\Gamma(\lceil\alpha\rceil - \alpha, imt_j)}{\Gamma(\lceil\alpha\rceil - \alpha)}\right)e^{imt_j - x_k^2}.
$$
Thus, we can discretize \eqref{e:PDE1} without its initial data as
$$
\Dta\cdot\U = \U\cdot\G^T + \HH,
$$
and, in order to incorporate the initial condition, we replace $\U$ by \eqref{e:Uboundaryt}:
$$
\Dta\cdot(\E\cdot\Uin+ \F) = (\E\cdot\Uin+ \F)\cdot\G^T + \HH.
$$
Then, left-multiplying by $\E^T$, and rearranging the terms, we get
$$
\E^T\cdot\Dta\cdot\E\cdot\Uin - \Uin\cdot\G^T =- \E^T\cdot\Dta\cdot\F + \E^T\cdot\HH,
$$
where we have used that $\E^T\cdot\E = \I$, and that $\E^T\cdot\F = \mathbf 0$. Defining $\A= \E^T\cdot\Dta\cdot\E$, $\B = -\G^T$, $\C = -\E^T\cdot\Dta\cdot\F + \E^T\cdot\HH$, we solve $\A\cdot\Uin + \Uin\cdot\B = \C$ by means of the Matlab command \verb|lyap|, and apply \eqref{e:Uboundaryt} to obtain $\U$. In Listing~\ref{code:pde2D}, we offer the Matlab implementation, corresponding to $\alpha = 0.97$ (i.e., a value of $\alpha$ close to 1), $m = 330$, $N_x = 16$, $N_t = 400$, $\dig = 400$ and $T = 2$. Note that, due to the highly oscilating nature of the solution $u(t, x)$ in time, it is obligatory to take a large enough matrix $\Dta$. The code takes 158.48 seconds to run, and $\mx{\U_{exact} - \U} = 6.1766\times10^{-13}$, which, in our opinion, is a pretty remarkable result.

\begin{lstlisting}[label=code:pde2D, language=Matlab, basicstyle=\footnotesize, caption = {Matlab program \texttt{pde2D.m}, that solves numerically }]
clear
tic
a=0.97; % $\alpha$
m=330; % $m$
Nx=16; % $N_x$
Nt=400; % $N_t$
dig=400; % Number of digits
T=2; % $T$
[~,Dta,t]=CaputoMatrix(Nt,a,T,dig);
b=1.4; % Scale factor $b$
[x,DD]=herdif(Nx,2,b);
Dx=DD(:,:,1); % $\mathbf D_x$
Dx2=DD(:,:,2); % $\mathbf D_x^2$
[tt,xx]=ndgrid(t,x); % Two-dimensional grid
u0=exp(-x.^2).'; % $u_0$
u=exp(1i*m*t-xx.^2); % $u$
E=[eye(Nt);zeros(1,Nt)]; % $\mathbf E$
F=[zeros(Nt,Nx);u0]; % $\mathbf F$
G=Dx2+2*x.*Dx+2*eye(Nx);
H=((1i*m)^a*(1-igamma(ceil(a)-a,1i*m*t)/gamma(ceil(a)-a))) ...
    .*exp(1i*m*tt-xx.^2); % $\mathbf H$
A=E.'*Dta*E; % $\mathbf A$
B=-G.'; % $\mathbf B$
C=-E.'*Dta*F+E.'*H; % $\mathbf C$
Uinner=lyap(A,B,-C); % $\mathbf U_{inner}$
U=E*Uinner+F; % $\mathbf U$
norm(u(:)-U(:),inf) % Error in max norm
toc % Elapsed time
\end{lstlisting}

\subsection{An example in multiple spatial dimensions}

The ideas in Section~\ref{s:oneD} can be extended to the case with multiple dimensions in space. To illustrate this, we consider the following Caputo-type advection-diffusion equation, which is a generalization of \eqref{e:PDE1}:
\begin{equation}
\label{e:PDEND}
	\left\{
	\begin{aligned}
		& \Dta u(t, \x) = \Delta u(t, \x) + 2\x\cdot\nabla u(t, \x) + 2Nu(t, \x)
			\cr		
				& \qquad\qquad\qquad + (im)^\alpha\left(1 - \frac{\Gamma(\lceil\alpha\rceil - \alpha, imt)}{\Gamma(\lceil\alpha\rceil - \alpha)}\right)e^{imt - \x\cdot\x}, & & (t,x)\in[0, 2]\times\mathbb R^d,
			\cr
		& u(0, \x) = u_0(\x) = e^{-\x\cdot\x},
	\end{aligned}
	\right.
\end{equation}
where $\alpha\in(0,1)$, $m\in\N$, the diffusive term is $\Delta u = \partial^2u/\partial x_1^2 + \partial^2u/\partial x_2^2 + \ldots + \partial^2u/\partial x_d^2$, and the advective term is $2\x\cdot\nabla u = 2(x_1\partial u/\partial x_1, x_2\partial u/\partial x_2, \ldots, x_d\partial u/\partial x_d)$. Moreover, the solution of \eqref{e:PDEND} is $u(t, \x) = e^{imt -\x\cdot\x}$. Then, similarly as in Section~\ref{s:oneD}, we discretize the operator $\cpt$ by means of $\Dta$, and the first and second spatial derivatives by means of $\D_{x_j}$ and $\D_{x_j}^2$, for $1\le j\le d$, where $\D_{x_j}\in\R^{N_j\times N_j}$ is a Hermite differentiation matrix generated by the Matlab function \verb|herdif|. This gives rise to an equation of the form
\begin{equation}
	\label{e:sylvesterND}
	\sum_{j = 1}^{d+1}\A_j\square_j\X = \A_1\square_1\X + \A_2\square_2\X + \ldots + \A_{d+1}\square_{d+1}\X = \B, \quad N \in \N,
\end{equation}
where $\A_j \equiv (a_{j,i_1i_2})\in\CC^{n_j\times n_j}$, $\B\equiv(b_{i_1i_2\ldots i_{d+1}})\in\CC^{n_1\times n_2\times\ldots\times n_{d+1}}$, $n_j\in\N$, for all $j$, and $\square_j$ indicates that the sum is performed along the $j$th dimension of $\X = (x_{i_1i_2\ldots i_{d+1}})\in \CC^{n_1\times n_2\times\ldots\times n_{d+1}}$, i.e.,
\begin{equation*} 
	[\A_j\square_j\X ]_{i_1i_2\ldots i_{N+1}} = \sum_{k = 1}^{n_j}a_{j,i_jk}x_{i_1i_2\ldots i_{j-1}k i_{j+1}\ldots i_{d+1}}.
\end{equation*}
Some references on the solution of Sylvester tensor equations like \eqref{e:sylvesterND} for three dimensions include \cite{benwen2010}, where the authors claimed to have solved such equation for the first time, and later on, \cite{delahozvadillo2013a,delahozvadillo2013b}. The references for higher dimensions are more recent, and we can mention, e.g., \cite{xinfang2021}, where a number of iterative algorithms based on Krylov spaces were applied; \cite{YuhanChen2023}, where a tensor multigrid method and an iterative tensor multigrid method was proposed; \cite{chen2020}, where recursive blocked algorithms were used; or \cite{massei2024}, where a a nested divide-and-conquer scheme was proposed.  However, in this section, we use the function \verb|sylvesterND| developed in \cite{cuestadelahoz2024}, that solves \eqref{e:sylvesterND} by means of an iterative, non-recursive method based on the Bartels-Stewart algorithm \cite{bartels1972}, that relies on the Schur decomposition of the matrices $\A_j$, and whose implementation is reduced to one single for-loop, and hence, independent from the number of dimensions.

Coming back to \eqref{e:PDEND}, let us denote $\U \equiv [u_{t_j,i_1,\ldots,i_d}]\in\CC^{(N_t+1)\times n_1\times\ldots\times n_d}$ the array such that $u_{t_j,i_1\ldots i_d} \approx u(t_j, x_{1,i_1}, \ldots, x_{d,i_d})$, and $\Uin\in\CC^{N_t\times n_1\times\ldots\times n_d}$ the array $\U$ without the entries of the form $u_{N_t+1,i_1,\ldots, i_d}$. Then, defining the matrix $\E$ as in \eqref{e:EF}, and the array $\F\equiv[f_{t_j,i_1,\ldots,i_d}]\in\CC^{(N_t+1)\times n_1\times\ldots\times n_d}$ consisting of zeros, except for the entries of the form $f_{N_t+1, i_1, \ldots, i_d}$, which are given by $f_{N_t+1, i_1, \ldots, i_d} = u_0(x_{i_1}, \ldots, x_{i_d})$, we have the following generalization of \eqref{e:Uboundaryt}:
\begin{equation}
\label{e:UboundarytND}
\U = \E\square_1\Uin+ \F.
\end{equation}
Moreover, denoting $\x_j = (x_{j,0}, x_{j,1}, \ldots, x_{j,n_j-1})$, and $\diag(\x_j)$ the diagonal matrix whose diagonal entries are precisely those of $\x_j$, for $1 \le j \le d$, we define
\begin{equation*}
\G_j = \D_{x_j}^2 + 2\diag(\x_j)\D_{x_j} + 2\I\in\R^{n_j\times n_j},
\end{equation*}
and $\HH \equiv [h_{t_j,i_1,\ldots,i_d}]\in\CC^{(N_t+1)\times n_1\times\ldots\times n_d}$, such that
$$
h_{t_j,i_1,\ldots,i_d} = (im)^\alpha\left(1 - \frac{\Gamma(\lceil\alpha\rceil - \alpha, imt_j)}{\Gamma(\lceil\alpha\rceil - \alpha)}\right)e^{imt_j - x_{i_1}^2 - \ldots x_{i_d}^2}.
$$
Thus, we can discretize \eqref{e:PDE1} without its initial data as
$$
\Dta\square_1\U = \sum_{j=1}^{d}\G_j\square_{j+1}\U + \HH,
$$
and, in order to incorporate the initial condition, we replace $\U$ by \eqref{e:UboundarytND}:
$$
\Dta\square_1(\E\square_1\Uin+ \F) = \sum_{j=1}^{d}\G_j\square_{j+1}(\E\square_1\Uin+ \F) + \HH.
$$
Then, left-multiplying by $\E^T$ along the first direction:
$$
\E^T\square_1(\Dta\square_1(\E\square_1\Uin+ \F)) = \sum_{j=1}^{d}\E^T\square_1(\G_j\square_{j+1}(\E\square_1\Uin+ \F)) + \E^T\square_1\HH,
$$
i.e.
$$
(\E^T\cdot\Dta\E)\square_1\Uin - \sum_{j=1}^{d}\G_j\square_{j+1}\Uin = -(\E^T\cdot\Dta)\square_1\F + \E^T\square_1\HH,
$$
where we have used that $\E^T\cdot\E = \I$, and that $\E^T\square_1\F = \mathbf 0$.
Defining $\A_1= \E^T\cdot\Dta\cdot\E$, $\A_{j+1} = -\G_j$, for $1\le j\le d$, and $\C = -(\E^T\cdot\Dta)\square_1\F + \E^T\square_1\HH$, we solve 
\begin{equation}
\sum_{j = 1}^{d+1}\A_j\square_j\Uin = \B,
\end{equation}
by means of the Matlab function \verb|sylvesterND| from \cite{cuestadelahoz2024}, and apply \eqref{e:UboundarytND} to obtain $\U$. Note that we have defined the Matlab function \verb|multND|, that computes $\square_1$, which is required to define $\B$, and to recover $\U$ from $\Uin$:
\begin{verbatim}
  function B=multND1(A,X)
  sizeX=size(X); % size of $\mathbf X$
  sizeB=[size(A,1),sizeX(2:end)]; % size of $\mathbf B$
  B=reshape(A*reshape(X,sizeX(1),[]),sizeB); % $\mathbf B$
\end{verbatim}
This function multiplies \verb|A| by the multidimensional array \verb|X|, which has been transformed into a two-dimensional matrix, and transforms the product back into a multidimensional array. The idea of transforming $\A$ into a two-dimensional matrix is taken from function \verb|multND| in \cite{cuestadelahoz2024}, which is invoked internally by \verb|sylvesterND| to multiply square matrices by multidimensional arrays, but, unlike \verb|multND1|, \verb|multND| requires to permute the dimensions of \verb|X| in a convenient way.

In Listing~\ref{code:pdeND}, we offer the Matlab implementation, corresponding to $\alpha = 0.97$, $m = 29$, $N_x = 16$, $N_t = 60$, $\dig = 60$, $T = 2$, and $d = 5$, i.e., 5 spatial dimensions, so the corresponding $\U$ is a six-dimensional array. We have programmed the code, so that it can be executed for any $d\in\N$. Note for instance the use of \verb|cell| structures to generate the $d+1$-dimensional grid:
\begin{verbatim}
  aux=[{t};repmat({x},d,1)];
  ttxx=cell(d+1,1);
  [ttxx{:}]=ndgrid(aux{:});
\end{verbatim}
or to define $\F$:
\begin{verbatim}
  F=zeros([Nt+1,Nx*ones(1,d)]); % $\mathbf F$
  index=[{Nt+1};repmat({':'},d,1)];
  F(index{:})=u0;
\end{verbatim}
The code takes 215.31 seconds to run, and $\mx{\U_{exact} - \U} = 6.5133\times10^{-14}$, which, again, we find it to be pretty remarkable.

\begin{lstlisting}[label=code:pdeND, language=Matlab, basicstyle=\footnotesize, caption = {Matlab program \texttt{pdeND.m}, that solves numerically }]
clear
tic
d=5; % Number of spatial dimensions
a=0.97; % $\alpha$
m=29; % $m$
Nx=16; % $N_x$
Nt=60; % $N_t$
dig=60; % Number of digits
T=2; % $T$
[~,Dta,t]=CaputoMatrix(Nt,a,T,dig);
b=1.4; % Scale factor $b$
[x,DD]=herdif(Nx,2,b);
Dx=DD(:,:,1); % $\mathbf D_x$
Dx2=DD(:,:,2); % $\mathbf D_x^2$
aux=[{t};repmat({x},d,1)];
ttxx=cell(d+1,1);
[ttxx{:}]=ndgrid(aux{:}); % (d+1)-dimensional grid
xx=cell(d,1);
aux=aux(2:end);
[xx{:}]=ndgrid(aux{:});
u0=1; % Initial data $u_0$
u=exp(1i*m*t); % Exact solution $u$
for j=1:d
    u0=u0.*exp(-xx{j}.^2);
    u=u.*exp(-ttxx{j+1}.^2);
end
E=[eye(Nt);zeros(1,Nt)]; % $\mathbf E$
F=zeros([Nt+1,Nx*ones(1,d)]); % $\mathbf F$
index=[{Nt+1};repmat({':'},d,1)];
F(index{:})=u0;
H=((1i*m)^a*(1-igamma(ceil(a)-a,1i*m*t)/gamma(ceil(a)-a))).*exp(1i*m*t); % $\mathbf H$
for j=1:d
    H=H.*exp(-ttxx{j+1}.^2);
end
B=-multND1(E.'*Dta,F)+multND1(E.',H); % $\mathbf B$
G=Dx2+2*x.*Dx+2*eye(Nx); % $\mathbf G$
AA=cell(d+1,1); % Store the matrices $\mathbf A_j$
AA{1}=E.'*Dta*E;
for j=1:d
    AA{j+1}=-G;
end
Uinner=sylvesterND(AA,B); % $\mathbf U_{inner}$
U=multND1(E,Uinner)+F; % $\mathbf U$
norm(U(:)-u(:),inf) % Error in max norm
toc % Elapsed time
\end{lstlisting}

\section*{Funding}

Francisco de la Hoz was partially supported by the research group grant IT1615-22 funded by the Basque Government, by the project PID2021-126813NB-I00 funded by MICIU/AEI/10.13039/501100011033 and by ``ERDF A way of making Europe'', and by the project PID2024-158099NB-I00, funded by MICIU. Peru Muniain was partially supported by the research group grant IT1461-22 funded by the Basque Government, and by the project PID2022-139458NB-I00 funded by MICIU.

\end{document}